\documentclass[hidelinks,onefignum,onetabnum]{siamart251216}

\usepackage{amsfonts,amssymb}
\usepackage{bm}
\usepackage{graphicx}
\usepackage{booktabs}
\usepackage{multirow}
\usepackage{makecell}
\usepackage{longtable}
\makeatletter
\def\LT@makecaption#1#2#3{\LT@mcol\LT@cols c{\hbox to\z@{\hss\parbox[t]\LTcapwidth{%
  \def\@captype{table}\@makecaption{#2}{#3}}\hss}}}
\makeatother
\usepackage{xcolor}
\usepackage{algorithmic}
\usepackage{tikz}
\usetikzlibrary{arrows.meta,positioning,calc,fit}
\usepackage{pgfplots}
\pgfplotsset{compat=1.18, table/comment chars={\#}}
\definecolor{siscgrey}{HTML}{A6A6A6}%
\definecolor{siscblue}{HTML}{1984C5}%
\definecolor{siscred}{HTML}{9E2A1F}%
\pgfplotsset{
  sisc/axis/.style={
    grid=major,
    major grid style={line width=0.4pt, draw=black!15},
    axis line style={line width=0.6pt, black},
    tick style={line width=0.6pt, black},
    tick label style={font=\small, black},
    label style={font=\small, black},
    title style={font=\small, black},
    legend cell align=left,
    legend style={font=\small, text=black, draw=none, fill=none,
      /tikz/every even column/.append style={column sep=0.9em}},
    legend image code/.code={\draw[mark repeat=2, mark phase=2]
      plot coordinates {(0cm,0cm) (0.42cm,0cm) (0.84cm,0cm)};},
    unbounded coords=jump,
    table/col sep=comma,
  },
  sisc/legend above/.style={
    legend style={at={(0.5,1.03)}, anchor=south},
  },
  sisc/ranks/.style={log basis x=2, xticklabel style={/pgf/number format/fixed}},
  sisc/seconds/.style={ytick={0.02,0.05,0.1,0.2,0.5,1,2,5},
    yticklabels={0.02,0.05,0.1,0.2,0.5,1,2,5}},
  sisc/series/.style={line width=1pt, mark size=2.4pt,
    mark options={solid, line width=0.6pt, draw=black, fill=#1}},
  sisc/jacobi/.style={sisc/series=siscgrey, color=siscgrey, solid, mark=*},
  sisc/cjpmg/.style={sisc/series=siscblue, color=siscblue, solid,
    mark=square*, mark size=2.1pt},
  sisc/hsfd/.style={sisc/series=siscred, color=siscred, dashed,
    mark=triangle*, mark size=2.9pt},
  sisc/fd/.style={sisc/series=white, color=black, dash dot,
    mark=diamond*, mark size=2.9pt},
}

\usetikzlibrary{arrows.meta,positioning,calc,fit,decorations.pathreplacing,patterns}
\definecolor{mfblue}{HTML}{1F77B4}
\definecolor{mforange}{HTML}{EB811B}
\definecolor{deepnavy}{HTML}{1C3A5E}

\newcommand{\siscdata}{figures/sisc_figures}

\ifpdf
  \DeclareGraphicsExtensions{.pdf,.png,.jpg}
\fi

\newsiamremark{remark}{Remark}

\newcommand{\bff}{\bm{f}}
\newcommand{\br}{\bm{r}}
\newcommand{\bu}{\bm{u}}

\newcommand{\bx}{\bm{x}}

\newcommand{\bs}{\bm{s}}

\newcommand{\bz}{\bm{z}}
\newcommand{\bzero}{\mathbf{0}}
\newcommand{\bone}{\mathbf{1}}
\newcommand{\bA}{\mathbf{A}}
\newcommand{\bB}{\mathbf{B}}
\newcommand{\bC}{\mathbf{C}}
\newcommand{\bD}{\mathbf{D}}
\newcommand{\bE}{\mathbf{E}}
\newcommand{\bI}{\mathbf{I}}
\newcommand{\bK}{\mathbf{K}}
\newcommand{\bM}{\mathbf{M}}
\newcommand{\bP}{\mathbf{P}}
\newcommand{\bR}{\mathbf{R}}
\newcommand{\bS}{\mathbf{S}}
\newcommand{\bW}{\mathbf{W}}

\newcommand{\bLam}{\bm{\Lambda}}
\newcommand{\order}{\mathcal{O}}
\DeclareMathOperator{\diag}{diag}

\headers{Matrix-free solvers for high-order finite elements}{G. Panigrahi and P. Motamarri}

\title{Accelerating matrix-free fast diagonalization and preconditioned linear solvers for
high-order finite element discretized elliptic problems on distributed GPUs\thanks{Submitted to the editors October 31, 2026.
\funding{G.P. was supported by the Prime Minister's Research Fellowship (PMRF), Ministry of
Education, India. P.M. was supported by the Google India Research Award 2023 and the Anusandhan
National Research Foundation (ANRF) ARG MATRICS grant (No. ANRF/ARGM/2025/002824/TS).}}}

\author{Gourab Panigrahi\thanks{Department of Computational and Data Sciences, Indian
  Institute of Science, Bengaluru 560012, India
  (\email{pgourab@iisc.ac.in}).}
\and Phani Motamarri\thanks{Department of Computational and Data Sciences, Indian
  Institute of Science, Bengaluru 560012, India
  (\email{phanim@iisc.ac.in}).}}

\newif\ifsiscappendix
\siscappendixfalse

\makeatletter
\newcommand*{\addFileDependency}[1]{\typeout{(#1)}\@addtofilelist{#1}\IfFileExists{#1}{}{\typeout{No file #1.}}}
\makeatother
\ifsiscappendix\else
\addFileDependency{supplement.tex}\addFileDependency{supplement.aux}
\fi

\ifsiscappendix
\usepackage{siscappendix}
\fi

\ifpdf
\hypersetup{
  pdftitle={Accelerating matrix-free fast diagonalization and preconditioned linear solvers for high-order finite element discretized elliptic problems on distributed GPUs},
  pdfauthor={G. Panigrahi and P. Motamarri}
}
\fi

\usepackage{pdfpages}

\begin{document}

\maketitle

\begin{abstract}
The time to solution of an iterative solver depends on the number of iterations and the cost per iteration. This paper develops algorithms that reduce both for elliptic problems discretized with high-order continuous finite elements on distributed GPUs, implemented in a unified matrix-free code supporting CUDA, HIP and SYCL. A distributed direct fast-diagonalization (FD) Poisson solver provides exact solutions on tensor-product meshes with direction-wise grading and coefficients varying along one coordinate direction. Its tensor data remain distributed, moved between slab layouts using collective transposes or kept in place by a Schur-complement approach. A two-level hybrid Schwarz preconditioner with FD cell solves (HS-FD) is extended to 2:1 hanging-node octrees with exact local cell inverses, while a \(p\)-multigrid preconditioner with Chebyshev--Jacobi smoothing (CJ-\(p\)MG) resolves hanging-node constraints at every level within fused kernels. Both approaches exploit FP32, where cell shapes and coefficients allow, within an FP64 conjugate-gradient iteration, adaptively select replicated or distributed coarse levels, support block solves for multiple right-hand sides and determine all parameters from the operator without case-specific tuning. Comparisons with Jacobi-preconditioned conjugate gradients (Jacobi PCG), MFEM, nekRS and deal.II span uniform, graded and octree meshes, layered coefficients, and four applications across NVIDIA, AMD and Intel GPUs. The direct solver is two to three orders of magnitude faster than Jacobi PCG on uniform meshes. HS-FD converges in 11 iterations for a generalized Poisson problem on octree meshes across 8 to 1024 MI250X GPUs. These results provide practical solver-selection guidelines based on mesh structure and coefficient characteristics. 
\end{abstract}

\begin{keywords}
matrix-free, $p$-multigrid, hybrid Schwarz, fast diagonalization, high-order finite elements, GPU
\end{keywords}

\begin{MSCcodes}
65F05, 65F08, 65F10, 65N30, 65N55, 65Y05, 65Y10, 65Y20
\end{MSCcodes}

\section{Introduction}
\label{sec:introduction}

Elliptic solves account for a large share of the run time in many simulation codes. For
instance, in incompressible flow an operator-splitting method solves a pressure Poisson problem at
every time step~\cite{Fischer1997Schwarz}, and in full-core reactor
simulations with the spectral-element code nekRS this pressure step typically takes
60--80\% of the run time~\cite{Min2022FullCore}. Implicit time stepping of viscous or
diffusive terms gives a screened Poisson operator,
the stiffness matrix plus a positive multiple of the mass
matrix~\cite{Min2022FullCore}. The generalized Poisson operator instead carries a variable
coefficient, which enters as the inverse density in two-phase flow~\cite{MacLachlan2008BubblyFlow}, the permeability of Darcy
flow~\cite{Bastian2019MatrixFreeDG}, the permittivity of implicit
solvation~\cite{Fisicaro2016GeneralizedPoisson} and the tissue conductivity of the forward
problem of electroencephalography~\cite{Vorwerk2014}. Capacitance extraction of interconnects solves one Poisson problem per
conductor~\cite{Nabors1991FastCap}, and resistivity tomography one conductivity problem per current
electrode~\cite{Jaysaval2023PFLOTRAN}. In Kohn--Sham density functional theory the
electrostatic potential is the solution of a Poisson problem in every self-consistent
iteration~\cite{Motamarri2020DFT-FECalculations}. In each of these applications the
same operator is solved many times on a fixed mesh, often for several right-hand sides.

The finite-element (FE) method handles complex geometries and locally refined meshes, and its
accuracy improves systematically as the mesh is refined and the polynomial degree is raised. This
work uses continuous ($C^0$) Lagrange elements of polynomial degree $p$ from 2 to 7 to discretize the underlying elliptic operators, and solves the
resulting sparse system by a preconditioned Krylov method, whose iterations are dominated by the
operator action and the preconditioner application. At high degree the operator action is
evaluated matrix-free by sum factorization, and fused kernels make it efficient on
GPUs~\cite{KronbichlerLjungkvist2019}.

Most preconditioners for sparse linear systems, including Gauss--Seidel relaxation, incomplete
factorizations and algebraic multigrid (AMG), are built from the entries of the assembled sparse
matrix. A matrix-free code never forms this matrix, so these preconditioners are not directly
available to it. AMG can still be applied to the sparse matrix of a simpler operator, such as a
low-order-refined one~\cite{Pazner2020LOR}, but its setup then adds to the cost. Assembling
the sparse matrix of the degree-$p$ discretization instead needs far more memory (about
100\,GB for $1.1\times10^{7}$ degrees of freedom (DoFs) at $p = 7$).
A matrix-free preconditioner is therefore built mainly from what the operator kernels already
compute: the operator action, the diagonal of the operator and the tensor-product structure of
the shape functions.

From the operator action, the diagonal and the tensor-product structure, we design three
solvers, one from each of three families, and order them by
how much of the inverse of the operator they apply: the diagonal, exact cell inverses or the exact inverse of the whole operator.
The $p$-multigrid preconditioner CJ-$p$MG applies the diagonal as a Chebyshev--Jacobi smoother on every level of a
hierarchy of polynomial degrees. Every level runs on the GPU in the sum-factorized kernel at a
lower degree on the cells and partition of the fine mesh, so an adaptively refined octree needs no mesh
hierarchy and no repartitioning. The price is a
trilinear coarse level on the fine mesh, replicated on every MPI rank except when BoomerAMG solves more than
$10^{6}$ unknowns or FD solves it on the Aurora 3D-IC stacks (\cref{sec:precond-common}).
The two-level hybrid Schwarz preconditioner HS-FD applies on every cell the exact inverse of a
cell problem by fast diagonalization (FD). Where the neighbors of a cell have its size, the cell
problem is the restriction of the operator to the cell. It remains a Kronecker sum of
one-dimensional operators next to hanging nodes and for a coefficient fitted by a separable
function (\cref{sec:asfd,sec:sm-precond}). The direct solve by FD applies the exact inverse of
the whole operator in the same way. It requires the whole operator to be a Kronecker sum, as for a coefficient that is
constant or varies along one direction, discretized on a tensor-product mesh, graded and fully
periodic ones included. Each method in this order inverts
more of the operator and needs more of its tensor structure. FD is the established local
solver of Schwarz methods for spectral elements because one FD solve costs about as much as one
operator action, so an iterative local solver would have to converge in a single iteration to
compete~\cite{Lottes2005HybridMultigrid}. This argument motivates FD here, in HS-FD and as
a direct solve on the whole mesh.

The three solvers are compared with Jacobi-preconditioned conjugate gradients (Jacobi PCG), the
solver of the CEED bake-off problems~\cite{Fischer2020BakeOff}, in the same code on
the same meshes, right-hand sides, tolerance and GPUs.
They are also compared with established solvers in three other codes: low-order-refined AMG in MFEM, the
$p$-multigrid of nekRS and the host $p$-multigrid of deal.II
(\cref{sec:results-amg,sec:results-iterations}).
Published speedups come from different codes, machines, tolerances and meshes,
so they do not tell which solver to prefer. This work asks which solver suits a distributed GPU code on a
given operator and mesh, and where its advantage comes from.

To find where the advantage comes from, the time to solution of each solver is split into the number of
iterations and the cost of one iteration. The number of iterations is set by the discrete
operator, the preconditioner setting, the Krylov method and the stopping test, so it is compared
across machines, and wall
times only within one machine (\cref{sec:results-setup}). The split shows where a
speedup comes from. For a generalized Poisson problem on an octree mesh (\cref{sec:results-iterations,fig:results-trp-scaling}) on 32 MI250X
graphics compute dies (GCDs), one HS-FD iteration costs as much as 14 Jacobi PCG iterations. HS-FD
needs 11 against 2317 to 2320 and is 15$\times$ faster.

\subsection{Related work}
\label{sec:related}

FD~\cite{Lynch1964TensorProduct} is used in high-order element
methods mostly as a local solver, in Schwarz methods and patch
smoothers~\cite{Lottes2005HybridMultigrid}, and as a global solver on
uniform meshes~\cite{Liu2023GPUSpectralPoisson} or for second-order finite
differences on graded grids. Costa et al.\ also apply it there on multiple GPUs~\cite{Costa2026NonUniformFD}.
The hybrid Schwarz--multigrid method of Lottes and Fischer~\cite{Lottes2005HybridMultigrid}
inverts by FD the operator restricted to an element extended by one layer of nodes, with counting
weights applied symmetrically so that the preconditioner can be used inside CG. nekRS runs
Chebyshev-accelerated smoothers of this kind on GPUs for conforming
meshes~\cite{Phillips2022TuningSEM}. Vertex patches give iterations nearly independent of the
polynomial degree, while cell patches with minimal overlap do
not~\cite{Brubeck2022VertexStar}, but FD inverts a vertex patch
only where its cells form a tensor-product grid, which excludes hanging nodes. The $p$-multigrid
method with Chebyshev--Jacobi smoothing is established for high-order
FE~\cite{Sundar2015MultigridComparison} and runs on hanging-node octrees on CPUs with an AMG coarse
solve~\cite{Munch2023}. High-order codes apply AMG only to the coarse problem of degree
1~\cite{Munch2023} or to a low-order-refined (LOR)
operator~\cite{Pazner2020LOR}. LOR with BoomerAMG runs on GPUs with hanging nodes, and for a
screened Poisson problem at $p = 6$ on one V100 GPU about 75\% of the time to solution is
setup~\cite{Pazner2023GPULOR} (\cref{sec:sm-related}).

The matrix-free $p$-multigrid and Schwarz preconditioners of these studies run on GPUs
only on conforming meshes~\cite{Fischer2022NekRS}, on equidistant Cartesian
meshes~\cite{Cui2025} or on a single GPU~\cite{KronbichlerLjungkvist2019}, and on hanging-node
octrees only on CPUs~\cite{Munch2023}. The global FD solves above are
restricted to uniform meshes, to one GPU or to second-order finite differences. In this work,
CJ-$p$MG and HS-FD run on hanging-node octrees on distributed GPUs. FD runs as a direct solve on
graded tensor-product meshes and as the exact coarse solve of both preconditioners. All three are
compared with Jacobi PCG on the same problems and GPUs, and with BoomerAMG on the
low-order-refined operator (LOR-AMG) on Polaris and Frontier.

\subsection{Contributions}
\label{sec:contributions}

The fused operator kernels, the GPU implementation of CG and the Chebyshev--Jacobi PCG used here
come from the GPU density functional theory framework of Kodali et al.~\cite{KodaliSC26} and are not
described again. An earlier GPU matrix-free Jacobi PCG of the authors was used for the
electrostatics of the 2023 Gordon Bell Prize study of Das et al.~\cite{DasGB2023}.
The use of shared memory and kernel fusion in the kernels of Kodali et al. follows the batched
matrix-free algorithms of Panigrahi et al.~\cite{PANIGRAHI2024104925}, and the block solves use a multilevel
batched layout similar to that of Panigrahi et al.~\cite{Panigrahi2026CPC}, adapted to GPUs. The
contributions are as follows.

\begin{itemize}
  \item A distributed-GPU implementation of CJ-$p$MG, the established $p$-multigrid with
    Chebyshev--Jacobi smoothing~\cite{Sundar2015MultigridComparison}, on hanging-node octrees, with the hanging-node constraints of every level resolved inside the
    fused operator kernels. The coarsest level is replicated on every MPI rank at one Allreduce
    per V-cycle, unless it is solved by BoomerAMG and has more than $10^{6}$ unknowns, when it is
    distributed. The fine and intermediate levels run in single precision except on
    strongly stretched cells or strong coefficient variation within a cell, and
    blocks of right-hand sides are solved together. On the octree generalized Poisson problem CJ-$p$MG needs 9 iterations on 8 to
    1024 Frontier GCDs and is about 9 to 13 times faster than Jacobi PCG.
  \item HS-FD, the weighted hybrid Schwarz--multigrid method of Lottes and
    Fischer~\cite{Lottes2005HybridMultigrid} with Chebyshev-accelerated FD cell
    solves~\cite{Phillips2022TuningSEM}, extended to 2:1 hanging-node octrees on distributed GPUs.
    Its weights vanish on the constrained nodes, so the preconditioner stays symmetric positive
    definite. A variable coefficient enters the cell problems through a separable fit on each cell, with a
    nodal scaling where the fit deviates from the coefficient. On strongly stretched cells the trilinear coarse level is
    replaced automatically by BoomerAMG on the low-order-refined operator. On the most strongly
    graded mesh, with cells of aspect ratio up to 16384, HS-FD then needs 18 to 23 iterations against 250
    for CJ-$p$MG. Its iterations grow with $p$.
  \item A direct FD solve on distributed GPUs. It inverts the FE operator exactly on
    tensor-product meshes with independent grading and Dirichlet or periodic conditions in each
    direction, also for a coefficient that varies along one direction. For the fully periodic
    operator it returns the solution of zero integral. The memory per MPI rank stays at
    $\order(N/P)$, right-hand sides are solved in blocks, and FD also serves as the exact
    coarse solve of both preconditioners on tensor-product meshes. On 36 Aurora tiles it is 37, 10 and
    $2.5\times10^{2}$ times faster than HS-FD on periodic uniform meshes and a graded mesh.
  \item A comparison of the three solvers with Jacobi PCG in one matrix-free code on the GPUs of
    three machines, NVIDIA A100 (Polaris), AMD MI250X (Frontier) and Intel Data Center GPU Max 1550
    (Aurora), with every preconditioner parameter chosen automatically for the operator at hand.
    BoomerAMG on the low-order-refined operator through MFEM and the $p$-multigrid of nekRS
    are compared in time, the host $p$-multigrid of deal.II in iterations. The tests cover
    uniform, graded and hanging-node meshes, many right-hand sides, and strong and weak scaling
    to 3072 MPI ranks with up to $1.57\times10^{9}$ DoFs, and give a solver selection by
    operator and mesh class (\cref{sec:results-selection}).
\end{itemize}

\Cref{sec:problem} states the problems, \cref{sec:fastdiag} the direct FD solve and
\cref{sec:precond} the two preconditioners and their GPU implementation,
\cref{sec:results} the cost model and measurements, and \cref{sec:conclusions} the conclusions
and limitations of each method.

\section{Problem formulation}
\label{sec:problem}

\subsection{Model operators and discretization}
\label{sec:problem-operators}

Three elliptic problems on a box $\Omega \subset \mathbb{R}^3$ are
considered: the Poisson equation, the generalized Poisson equation with a
variable coefficient, and the screened Poisson equation,
\begin{equation}
  -\nabla^2 u(\bx) = g(\bx), \qquad
  -\nabla \cdot \bigl(\epsilon(\bx) \nabla u(\bx)\bigr) = g(\bx), \qquad
  -\nabla^2 u(\bx) + \sigma u(\bx) = g(\bx),
  \label{eq:prob-strong}
\end{equation}
where $u(\bx)$ is the unknown, $g(\bx)$ is the source, $\epsilon(\bx) \ge \epsilon_{\min} > 0$
is a scalar coefficient and $\sigma > 0$ is a constant. The screened Poisson equation arises when
diffusion terms are integrated implicitly in time, with $\sigma$ inversely proportional to the
time step.

The coefficient $\epsilon(\bx)$ comes from the application. In the test problems of
\cref{sec:results} it is the dielectric function of implicit
solvation~\cite{Fisicaro2016GeneralizedPoisson}, the inverse density of bubbly
flow~\cite{MacLachlan2008BubblyFlow}, the conductivity of a fiber
composite~\cite{Moulinec1998FFTHomogenization} or the tissue conductivity of a human
head~\cite{HuangNYHead}. In several problems the coefficient varies along one direction only,
$\epsilon(\bx) = \kappa(z)$, with $\kappa$ constant on layers. Such a layered coefficient
describes the earth in resistivity tomography~\cite{Jaysaval2023PFLOTRAN}, a layered Darcy medium
or the stack of a three-dimensional integrated circuit (3D-IC)~\cite{Guo2026ThermalLabyrinth}.

The domain is meshed with hexahedral cells, and $u$ is approximated by
continuous finite elements of polynomial degree $p$ in each direction (FE
degree $p$), with shape functions $\phi_J(\bx)$ on the $n_p = p + 1$
Gauss--Lobatto--Legendre (GLL) nodes per direction of each cell. The
integrals are evaluated with the tensor-product Gauss rule of
$n_q = n_p$ points per direction, so the quadrature points do not coincide
with the GLL nodes and the mass matrix is consistent (non-diagonal). The generalized Poisson
equation is discretized in its divergence form, with $\epsilon$ evaluated at the quadrature
points. Munch et al.\ use the same rule, with $p + 1$ Gauss points per direction, in
their matrix-free multigrid~\cite{Munch2023}. The CEED bake-off problems use $p + 2$
Gauss points, or the $p + 1$ GLL nodes themselves as in spectral-element codes~\cite{Fischer2020BakeOff}. On the box-shaped cells used here,
the rule computes the stiffness and mass matrices exactly where the coefficient is constant on
each cell. Where $\epsilon$ varies or jumps inside a cell, it is sampled at the quadrature points. This limits the accuracy of the discretization, not the comparison of
the solvers, which all solve the same discrete system.
The discrete problem is
\begin{equation}
  \bA \bu = \bff, \qquad
  \bA \in \{\bK,\ \bK^{\epsilon},\ \bK + \sigma \bM\},
  \label{eq:prob-system}
\end{equation}
where $\bu$ holds the nodal values of $u^h(\bx) = \sum_J u_J \phi_J(\bx)$, the
right-hand side is $f_J = \int_\Omega \phi_J(\bx)\, g(\bx)\, d\bx$, and
$K_{IJ} = \int_\Omega \nabla \phi_I \cdot \nabla \phi_J\, d\bx$,
$K^{\epsilon}_{IJ} = \int_\Omega \epsilon(\bx)\, \nabla \phi_I \cdot \nabla \phi_J\, d\bx$ and
$M_{IJ} = \int_\Omega \phi_I \phi_J\, d\bx$ are the stiffness, generalized stiffness
and mass matrices. The number of nodes of the mesh, constrained nodes included, is denoted
by $N$, and the number of MPI ranks by $P$. Each MPI rank drives one GPU, or one die of a GPU
with two dies (\cref{sec:results-setup}). The matrix $\bA$ is never assembled: its action on a
vector is evaluated cell by cell by sum factorization (\cref{sec:implementation}).

The simulation domain is discretized with two mesh classes. A tensor mesh has nodes that form the Cartesian
product of three one-dimensional node sets, each uniform or graded
independently, as for a box or a slab stretched along its normal. An octree
mesh is refined locally, with neighboring cells differing by at most one
level (2:1), and carries hanging nodes on the faces and edges between
refinement levels.

\subsection{Constraints and the singular periodic case}
\label{sec:problem-constraints}

A hanging node is a weighted sum of its constraining nodes on the coarser neighbor, a periodic
constraint copies the value from the constraining to the constrained node, and a homogeneous
Dirichlet node is zero. One constraint matrix $\bC$ maps the free DoFs $\hat{\bu}$ (the
unconstrained nodes, constraining nodes included) to all nodes, $\bu = \bC \hat{\bu}$, and testing
\cref{eq:prob-system} with the free DoFs gives the system that is solved,
\begin{equation}
  \bC^T \bA \bC\, \hat{\bu} = \bC^T \bff.
  \label{eq:prob-constraint}
\end{equation}
$\bC$ and its transpose are applied cell-locally inside the operator action, and from here on
$\bA$, $\bK$ and $\bM$ denote $\bC^T \bA \bC$, $\bC^T \bK \bC$ and $\bC^T \bM \bC$ on the free DoFs.

With periodic conditions in all three directions and $\sigma = 0$, the matrix $\bK$ is singular,
and its null space is spanned by the vector of ones, $\bone$. The system then has a solution
only if $\bone^T \bff = 0$, so the mean of the entries of $\bff$ is subtracted from $\bff$. The
solution is unique up to an added constant, which each method fixes in its own way. The direct
solve of \cref{sec:fastdiag} applies a pseudoinverse that returns the zero-integral solution,
$\int_\Omega u^h\, d\bx = \bone^T \bM \bu = 0$. Inside CG the constant is removed from the preconditioned
residual by a projection weighted by the diagonal of the operator, the metric-matched projection.
This weight matches the Jacobi and Chebyshev--Jacobi preconditioners, and for CJ-$p$MG and HS-FD
it leaves the iterations unchanged. Since $\bone^T\bff = 0$ and $\bone^T\bK = \bzero^T$, every
residual of CG is orthogonal to $\bone$, and a projection of any weight changes the preconditioned
residual only by a multiple of $\bone$. CG reads that vector only through its product with $\bK$,
which annihilates $\bone$, and its inner products with residuals, so the residuals and step
lengths are unchanged and the iterate moves only by a constant (\cref{sec:sm-metric}).
On the fully periodic generalized Poisson problems CG instead solves a system with a rank-one shift that
removes the null space, and the preconditioner is matched to the shift (\cref{sec:sm-periodic}).

\subsection{Iterative solver strategy}
\label{sec:problem-krylov}

All preconditioners are applied inside the conjugate gradient method (CG) on the GPU
implementation of Kodali et al.~\cite{KodaliSC26}, with the single-reduction recurrence of
Chronopoulos and Gear~\cite{ChronopoulosGear1989}. CJ-$p$MG and HS-FD run partly in FP32 and vary
slightly between applications, so for them CG uses the flexible (Polak--Ribi\`ere)
coefficient~\cite{Notay2000FlexibleCG}. Jacobi PCG uses the standard coefficient, and two-reduction
CG on the systems marked in the tables, where the single-reduction recurrence stalls. For a block of
right-hand sides, CJ-$p$MG and HS-FD solve all right-hand sides at once: each keeps its own step
lengths and search directions, the preconditioner is applied to all of them together, and their
inner products share the reductions of the iteration. Jacobi PCG on the generalized Poisson blocks of \cref{sec:results-mrhs} uses the
breakdown-free block CG of Ji and Li~\cite{JiLi2017BFBCG} (\cref{sec:sm-stopping}).

A solve stops when the relative true residual $\|\bff - \bA \bu\|_2 / \|\bff\|_2$ falls below
$10^{-10}$, or after at most $20\,000$ iterations. Capped solves, the backward error of
solves with a large coefficient ratio (the largest over the smallest coefficient value) and the two
problems solved to $10^{-9}$ are treated in \cref{sec:results-setup}.

\section{Direct solve by fast diagonalization}
\label{sec:fastdiag}

On a tensor mesh with a constant coefficient, the discrete operator of the Poisson and screened
Poisson equations \cref{eq:prob-strong} is a sum of Kronecker products of one-dimensional stiffness
and mass matrices, and fast diagonalization (FD) applies its inverse
exactly~\cite{Lynch1964TensorProduct}. For high-order elements FD has been used mostly as a local
solver, and globally only on uniform meshes, on one GPU or for finite differences
(\cref{sec:related}). The Kronecker structure, the exact inverse, the treatment of the null mode
and the transposes of the distributed solve follow known
constructions~\cite{Lynch1964TensorProduct,Costa2026NonUniformFD}. This work adds the direct FD
solve of the high-order FE operator on graded meshes on distributed GPUs, with $\order(N/P)$ memory
per MPI rank and blocks of right-hand sides, extends it to a coefficient that varies along one
direction (\cref{sec:fd-kz}), and uses it as the exact coarse solve of both multigrid
preconditioners (\cref{sec:fd-coarse}).

\subsection{Kronecker-sum structure}
\label{sec:fd-kron}

On a tensor mesh, every global shape function is a product of
one-dimensional FE shape functions, $\phi_{ijk} = \phi^x_i \phi^y_j \phi^z_k$, and the Gauss quadrature rule of \cref{sec:problem-operators} is a tensor
product as well. Every integral arising in \cref{eq:prob-system} then factorizes
into one-dimensional integrals, and the operator $\bA = \bK + \sigma \bM$,
with $\sigma \ge 0$, takes the form
\begin{equation}
  \bK = \bK_x \otimes \bM_y \otimes \bM_z + \bM_x \otimes \bK_y \otimes \bM_z
      + \bM_x \otimes \bM_y \otimes \bK_z, \qquad
  \bM = \bM_x \otimes \bM_y \otimes \bM_z,
  \label{eq:fd-kron}
\end{equation}
where $\bK_d, \bM_d \in \mathbb{R}^{n_d \times n_d}$ are the one-dimensional
stiffness and mass matrices along direction $d \in \{x, y, z\}$ and $n_d$ is
the number of free nodes along $d$. They are assembled over the $n^d_e$ cells along $d$, each
with its own size $h_e$, so grading a direction changes only the pair $(\bK_d, \bM_d)$. The
mass matrices are consistent, so FD diagonalizes each pair by a generalized eigenproblem
$\bK_d\bs = \lambda\bM_d\bs$ (\cref{sec:fd-inverse}) and inverts the discrete operator itself, with
no mass lumping, to rounding error.

The structure survives the constraints of \cref{eq:prob-constraint}
whenever $\bC = \bC_x \otimes \bC_y \otimes \bC_z$, as for homogeneous Dirichlet and periodic
conditions, which act on one direction at a time. A Dirichlet direction removes its two end nodes, so $n_d = n^d_e p - 1$ and
$\bK_d$ is symmetric positive definite. A periodic direction identifies its
last node with its first, so $n_d = n^d_e p$, the pair $(\bK_d, \bM_d)$
gains entries in its top-right and bottom-left corners, and $\bK_d$ has a
simple zero eigenvalue with the constant vector as eigenvector.
At hanging nodes a line of nodes crossing a refinement interface has more free nodes on the fine
side than on the coarse side, so the free DoFs are not a Cartesian product, and FD does not apply
to the whole operator.

Robin conditions with a positive coefficient constant on whole faces add it to the end diagonal
entry of $\bK_d$ and keep the structure (\cref{sec:sm-robin}), but the experiments use Dirichlet and
periodic directions.

The generalized Poisson operator keeps the Kronecker structure when its coefficient varies along
$z$ only and is constant on every cell layer, as in layered media such as the 3D-IC stack: the pair
$(\bK_z, \bM_z)$ is then replaced by one weighted by the coefficient of each layer
(\cref{sec:fd-kz}). A coefficient that does not separate into one-dimensional factors, such as the
coefficient of the octree problems of \cref{sec:results}, breaks the Kronecker structure, and the problem is
then solved by CG.

\subsection{Exact inverse}
\label{sec:fd-inverse}

For each direction, the one-dimensional generalized eigenproblem
\begin{equation}
  \bK_d \bS_d = \bM_d \bS_d \bLam_d, \qquad \bS_d^T \bM_d \bS_d = \bI,
  \label{eq:fd-evp}
\end{equation}
is solved once per mesh, where $\bLam_d = \diag(\lambda^d_0, \ldots,
\lambda^d_{n_d - 1})$ holds the eigenvalues in ascending order and the
columns of $\bS_d$ are $\bM_d$-orthonormal eigenvectors. With
$\bS = \bS_x \otimes \bS_y \otimes \bS_z$, the mixed-product rule applied to
\cref{eq:fd-kron} gives $\bS^T \bM \bS = \bI$  and further
\begin{equation}
  \bS^T \bA \bS = \bLam, \qquad
  \Lambda_{ijk} = \lambda^x_i + \lambda^y_j + \lambda^z_k + \sigma,
  \label{eq:fd-diag}
\end{equation}
where $\bLam$ is diagonal and $\Lambda_{ijk}$ is its entry for the index
triple $(i, j, k)$. If $\bA$ is nonsingular, its inverse is therefore
\begin{equation}
  \bA^{-1} = \bS \bLam^{-1} \bS^T
  = (\bS_x \otimes \bS_y \otimes \bS_z)\, \bLam^{-1}\,
    (\bS_x^T \otimes \bS_y^T \otimes \bS_z^T).
  \label{eq:fd-inverse}
\end{equation}
Every $\lambda^d_i$ is positive on a Dirichlet direction and non-negative on
a periodic one. So $\bA$ is singular only when all three directions are
periodic and $\sigma = 0$. Then exactly one diagonal entry vanishes,
$\Lambda_{000} = \lambda^x_0 + \lambda^y_0 + \lambda^z_0 = 0$, and its eigenvector is
constant.

For the fully periodic operator with $\sigma = 0$, the FD solve sets the reciprocal of the zero
eigenvalue to zero,
\begin{equation}
  \bA^{+}_{\bM} = \bS \bLam^{+} \bS^T, \qquad
  \Lambda^{+}_{ijk} =
  \begin{cases}
    0, & (i, j, k) = (0, 0, 0), \\
    1/\Lambda_{ijk}, & \text{otherwise},
  \end{cases}
  \label{eq:fd-nullmode}
\end{equation}
where $\bLam^{+}$ is the diagonal matrix of the entries $\Lambda^{+}_{ijk}$. Because $\bS$ is
$\bM$-orthonormal, $\bA^{+}_{\bM}$ is the pseudoinverse in the mass-matrix inner product, and it
returns the zero-integral solution of \cref{sec:problem-constraints}
(\cref{sec:fd-periodic}).

\begin{proposition}
\label{prop:fd-pinv}
Let all three directions be periodic, $\sigma = 0$, and
$\bu = \bA^{+}_{\bM} \bff$ as in \cref{eq:fd-nullmode}. Then, for every
$\bff$,
\begin{enumerate}
  \item $\int_\Omega u^h = \bone^T \bM \bu = 0$, and
  \item $\bA \bu = \bff - \dfrac{\bone^T \bff}{\bone^T \bM \bone}\, \bM \bone$,
    so $\bA \bu = \bff$ whenever $\bone^T \bff = 0$.
\end{enumerate}
\end{proposition}

\begin{proof}
The column of $\bS$ for the triple $(0,0,0)$ is the $\bM$-normalized
constant vector $\bS \bm{e}_0 = \bone / \sqrt{\bone^T \bM \bone}$, where
$\bm{e}_0$ is the corresponding unit vector. Since $\Lambda^{+}_{000} = 0$,
$\bu$ is a combination of the other columns of $\bS$, each of which is
$\bM$-orthogonal to $\bone$ by $\bS^T \bM \bS = \bI$, so $\bone^T \bM \bu =
0$. The partition of unity $\sum_I \phi_I \equiv 1$ gives $\bone^T \bM \bu =
\int_\Omega (\sum_I \phi_I)(\sum_J u_J \phi_J) = \int_\Omega u^h$, which proves
the first item. For the second, $\bS^T \bM \bS = \bI$ implies
$\bM \bS \bS^T = \bI$, and \cref{eq:fd-diag} implies $\bA \bS = \bM \bS
\bLam$. Hence $\bA \bu = \bM \bS \bLam \bLam^{+} \bS^T \bff = \bM \bS (\bI -
\bm{e}_0 \bm{e}_0^T) \bS^T \bff = \bff - \bM \bS \bm{e}_0 \bm{e}_0^T \bS^T \bff$, which
is the stated expression.
\end{proof}

\subsection{Algorithm and complexity}
\label{sec:fd-algorithm}

\begin{figure}[tbhp]
\centering
\def\fdbox#1#2#3#4#5#6#7{%
  \fill[#7] (#1,#2,#3)--(#4,#2,#3)--(#4,#2,#6)--(#1,#2,#6)--cycle;%
  \fill[#7] (#1,#2,#6)--(#4,#2,#6)--(#4,#5,#6)--(#1,#5,#6)--cycle;%
  \fill[white,opacity=0.35] (#1,#2,#6)--(#4,#2,#6)--(#4,#5,#6)--(#1,#5,#6)--cycle;%
  \fill[#7] (#4,#2,#3)--(#4,#5,#3)--(#4,#5,#6)--(#4,#2,#6)--cycle;%
  \fill[black,opacity=0.12] (#4,#2,#3)--(#4,#5,#3)--(#4,#5,#6)--(#4,#2,#6)--cycle;%
  \draw[deepnavy!45,line width=0.4pt,line join=round]%
    (#1,#2,#3)--(#4,#2,#3)--(#4,#2,#6)--(#1,#2,#6)--cycle%
    (#1,#2,#6)--(#1,#5,#6)--(#4,#5,#6)--(#4,#2,#6)%
    (#4,#2,#3)--(#4,#5,#3)--(#4,#5,#6);}%
\def\fdgridcube#1{%
  \fdbox{0}{0}{0}{#1}{#1}{#1}{mfblue!35}%
  \foreach \t in {0.25,0.5,0.75}{%
    \draw[white,line width=0.5pt] (\t*#1,0,0)--(\t*#1,0,#1)--(\t*#1,#1,#1);%
    \draw[white,line width=0.5pt] (0,0,\t*#1)--(#1,0,\t*#1)--(#1,#1,\t*#1);%
    \draw[white,line width=0.5pt] (0,\t*#1,#1)--(#1,\t*#1,#1)--(#1,\t*#1,0);}%
  \draw[deepnavy!45,line width=0.4pt] (0,0,0)--(#1,0,0)--(#1,0,#1)--(0,0,#1)--cycle%
    (0,0,#1)--(0,#1,#1)--(#1,#1,#1)--(#1,0,#1) (#1,0,0)--(#1,#1,0)--(#1,#1,#1);}%
\def\fdlinesx#1{\foreach \t in {0.25,0.5,0.75}{\draw[deepnavy,line width=1pt] (0,0,\t*#1)--(#1,0,\t*#1);}}%
\def\fdlinesy#1{\foreach \t in {0.25,0.5,0.75}{\draw[deepnavy,line width=1pt] (\t*#1,0,#1)--(\t*#1,#1,#1);}}%
\def\fdlinesz#1{\foreach \t in {0.25,0.5,0.75}{\draw[deepnavy,line width=1pt] (\t*#1,0,0)--(\t*#1,0,#1);}}%
\def\fdmini#1#2{\begin{scope}[shift={(#1,0)},fdthree]\fdgridcube{0.6}\csname fdlines#2\endcsname{0.6}\end{scope}}%
\begin{tikzpicture}[font=\footnotesize,every node/.append style={text=black},
    fdthree/.style={x={(1cm,0cm)}, y={(0.38cm,0.26cm)}, z={(0cm,1cm)}},
    fdflow/.style={-{Stealth[length=4pt]},deepnavy,line width=1pt},
    fdstep/.style={-{Stealth[length=3pt]},deepnavy,line width=0.6pt},
    fdbrace/.style={decorate,decoration={brace,amplitude=3pt},deepnavy,line width=0.6pt}]
  \node[anchor=west] at (-0.75,1.75) {(a)};
  \begin{scope}[shift={(-0.7,0.05)},fdthree]
    \draw[fdstep] (0,0,0)--(0.35,0,0) node[right,inner sep=1pt] {$x$};
    \draw[fdstep] (0,0,0)--(0,0.55,0) node[above right,inner sep=0pt] {$y$};
    \draw[fdstep] (0,0,0)--(0,0,0.35) node[above,inner sep=1pt] {$z$};
  \end{scope}
  \begin{scope}[fdthree]\fdgridcube{0.85}\end{scope}
  \node at (0.58,1.32) {$\mathcal{R}$, from $\bff$};
  \node at (0.58,-0.25) {$n_x\times n_y\times n_z$};
  \draw[fdflow] (1.25,0.45)--(1.62,0.45);
  \fdmini{1.7}{x}\fdmini{2.75}{y}\fdmini{3.8}{z}
  \draw[fdstep] (2.56,0.3)--(2.72,0.3); \draw[fdstep] (3.61,0.3)--(3.77,0.3);
  \node at (2.0,-0.22) {$\times_x \bS_x^T$};
  \node at (3.05,-0.22) {$\times_y \bS_y^T$};
  \node at (4.1,-0.22) {$\times_z \bS_z^T$};
  \draw[fdbrace] (1.7,0.95)--(4.63,0.95) node[midway,above=3pt] {mode products};
  \draw[fdflow] (4.7,0.45)--(5.05,0.45);
  \node[draw=deepnavy,line width=1pt,fill=mforange,rounded corners=2pt,minimum width=1.25cm,
        minimum height=0.75cm,align=center,inner sep=2pt] (scale) at (5.75,0.45)
        {scale by\\ $\Lambda^{-1}_{ijk}$};
  \node[align=center] at (5.75,1.35) {$\Lambda_{ijk} = \lambda^x_i + \lambda^y_j$\\ $+\, \lambda^z_k + \sigma$};
  \draw[fdflow] (6.45,0.45)--(6.82,0.45);
  \fdmini{6.9}{x}\fdmini{7.95}{y}\fdmini{9.0}{z}
  \draw[fdstep] (7.76,0.3)--(7.92,0.3); \draw[fdstep] (8.81,0.3)--(8.97,0.3);
  \node at (7.2,-0.22) {$\times_x \bS_x$};
  \node at (8.25,-0.22) {$\times_y \bS_y$};
  \node at (9.3,-0.22) {$\times_z \bS_z$};
  \draw[fdbrace] (6.9,0.95)--(9.83,0.95) node[midway,above=3pt] {mode products};
  \draw[fdflow] (9.9,0.45)--(10.27,0.45);
  \begin{scope}[shift={(10.35,0)},fdthree]\fdgridcube{0.85}\end{scope}
  \node at (10.93,1.32) {$\bu$};
  \node at (10.93,-0.25) {solution};
  \begin{scope}[shift={(0,-2.85)}]
    \node[anchor=west] at (-0.75,1.85) {(b)};
    \node[align=center] at (0.15,0.55) {locally\\ owned\\ DoFs};
    \begin{scope}[shift={(1.6,0)},fdthree]
      \foreach \r/\s in {0/15,1/35,2/55,3/75}{%
        \fdbox{0}{0}{0.2375*\r}{0.95}{0.95}{0.2375*\r+0.2375}{mfblue!\s}}
      \fdlinesx{0.95}
    \end{scope}
    \node at (2.25,1.55) {$k$-slab};
    \node[align=center] at (2.25,-0.42) {$\times_x$, $\times_y$ local};
    \draw[fdflow] (0.72,0.75)--(1.53,0.75) node[midway,above] {entry};
    \draw[fdflow] (1.53,0.3)--(0.72,0.3) node[midway,below] {exit};
    \draw[fdflow] (3.02,0.75)--(4.22,0.75) node[midway,above] {transpose};
    \draw[fdflow] (4.22,0.3)--(3.02,0.3) node[midway,below] {transpose};
    \begin{scope}[shift={(4.3,0)},fdthree]
      \foreach \r/\s in {0/15,1/35,2/55,3/75}{%
        \fdbox{0.2375*\r}{0}{0}{0.2375*\r+0.2375}{0.95}{0.95}{mfblue!\s}}
      \fdlinesz{0.95}
    \end{scope}
    \node at (4.95,1.55) {$i$-slab};
    \node[align=center] at (4.95,-0.5) {$\times_z \bS_z^T$, scaling,\\ $\times_z \bS_z$};
    \draw[deepnavy,line width=0.4pt,dashed] (6.35,-0.85)--(6.35,1.85);
    \begin{scope}[shift={(6.85,0)},fdthree]
      \foreach \iy/\ix/\s in {1/0/55,1/1/75,0/0/15,0/1/35}{%
        \fdbox{0.475*\ix}{0.475*\iy}{0}{0.475*\ix+0.475}{0.475*\iy+0.475}{0.95}{mfblue!\s}}
      \fdlinesz{0.95}
    \end{scope}
    \node at (7.5,1.55) {rank grid};
    \node[align=center] at (7.5,-0.5) {lines along $z$\\ local};
    \node[align=left,anchor=west] at (8.45,0.55) {ranks outnumber\\ the cells along $z$:\\ two-dimensional\\ grid of ranks,\\ transposes within\\ one row or column};
  \end{scope}
\end{tikzpicture}
\caption[Direct FD solve]{Direct FD solve. (a) The right-hand side tensor $\mathcal{R}$ is
multiplied along its lines by $\bS_x^T$, $\bS_y^T$ and $\bS_z^T$, scaled entrywise by
$\Lambda^{-1}_{ijk}$ and multiplied by $\bS_x$, $\bS_y$ and $\bS_z$. (b) The distributed solve
applies the products along $x$ and $y$ in slabs along $z$ ($k$-slab) and those along $z$ in slabs
along $x$ ($i$-slab), with a transpose between the two layouts, or on a two-dimensional grid of
MPI ranks when the ranks outnumber the cells along $z$. Shading marks four MPI ranks.}
\label{fig:fd-solve}
\end{figure}
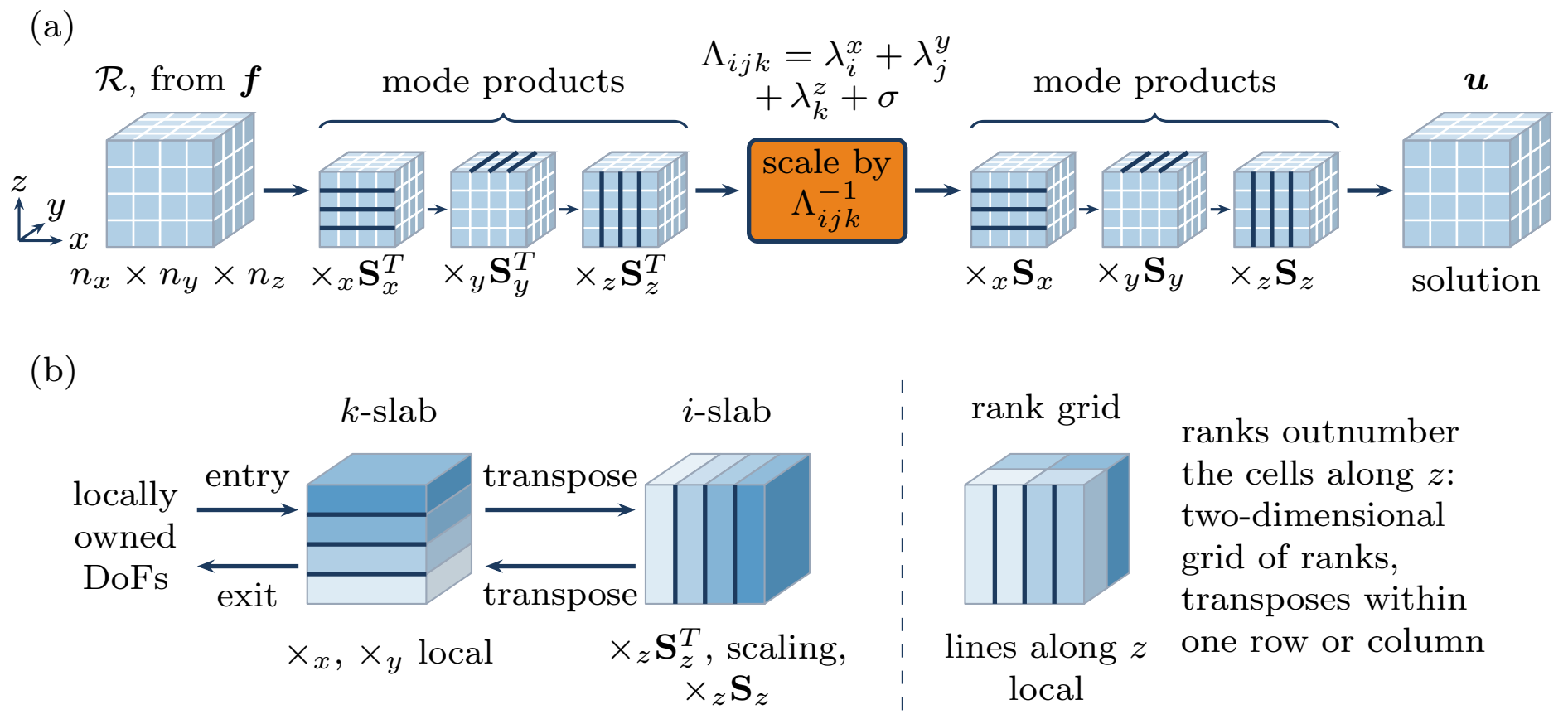

The free entries of $\bff$ are arranged as a tensor
$\mathcal{R} \in \mathbb{R}^{n_x \times n_y \times n_z}$. A line of the
tensor along $x$ holds the values on one row of nodes parallel to the $x$ axis (the entries
$R_{ijk}$ with $j$ and $k$ fixed), and lines along $y$ and $z$ are defined in the same way. The Kronecker factors of
\cref{eq:fd-inverse} are applied as mode
products~\cite{Kolda2009TensorDecompositions}. The mode-$d$ product
$\times_d$ multiplies every line along direction $d$ by an $n_d \times n_d$
matrix, which is one dense matrix--matrix product with the unfolding whose columns are these lines. A solve applies the three mode products by $\bS_d^T$, scales every entry by
$1/\Lambda_{ijk}$ and applies the three products by $\bS_d$ (\cref{fig:fd-solve}(a),
\cref{alg:fd-solve}).

One solve on $n$ nodes per direction costs $12 n^4 + n^3 = \order(N^{4/3})$ operations in six
dense matrix--matrix products and one pointwise scaling (\cref{sec:sm-schur}). The setup is done
once per mesh, before the first solve. It checks that the mesh is a tensor product, to within
$10^{-4} h_e$ in the node coordinates, and that the free DoFs form a Cartesian product. For the
generalized Poisson operator it also checks that the coefficient is constant on every cell layer
along $z$ (\cref{sec:fd-kz}). A problem that fails any of these tests is solved by CG. Otherwise
the setup solves the three eigenproblems of \cref{eq:fd-evp}, at $\order(n^3)$ operations each.
The solve runs entirely in double precision.

\subsection{Distributed solve}
\label{sec:fd-distributed}

This work proposes the following strategy for the FD solve on many MPI ranks
(\cref{fig:fd-solve}(b)). A mode product needs complete lines of the tensor along its direction.
The tensor is therefore held in one of two layouts and moved between them: slabs along $z$
($k$-slabs), in which every rank holds complete lines along $x$ and $y$, and slabs along $x$
($i$-slabs), with complete lines along $z$. Each rank stores $\order(N/P)$ entries of the tensor
and a copy of the one-dimensional eigenvectors. The partition of cells and DoFs used by the FE
code is not changed.

In the general case a solve takes four collective steps. The entry map moves the locally owned
DoFs into $k$-slabs, where the products along $x$ and $y$ are applied. A transpose moves the
tensor to $i$-slabs for the products along $z$ and the scaling. A second transpose and the exit
map return the result to the locally owned DoFs. The solver of Costa et al.\
uses such collective transposes on GPUs~\cite{Costa2026NonUniformFD}. If the FE cell partition
matches the $k$-slabs, the entry and exit maps reduce to local copies.

After the products along $x$ and $y$, the solve along $z$ splits into one banded system per line
along $z$. When $z$ can be cut at cell boundaries with at least one cell per rank, two on a
periodic direction, each rank solves these systems on its own cells, and the values on the nodes
at the cuts follow from a tridiagonal Schur complement system. This replaces both transposes, and
each exchange moves $2 n_x n_y$ values per rank in place of $n_x n_y n_z / P$. For a block of $n_v$ right-hand sides, 
each mode product becomes one matrix product with $n_v$ times as many columns, and each exchange
carries all $n_v$ vectors in one message. A replicated
scheme, which all-gathers the right-hand side and solves redundantly on every rank, is slower
than the distributed one on two or more MPI ranks (\cref{sec:results}). \Cref{sec:sm-schur}
gives the partition and the Schur complement solve.

FD also gives the exact coarse solve of CJ-$p$MG and HS-FD wherever their trilinear coarse grid is
a tensor grid and the coefficient varies at most along $z$ (\cref{sec:fd-coarse}).

\section{Preconditioners}
\label{sec:precond}

This section describes the two preconditioners: CJ-$p$MG, the Chebyshev--Jacobi $p$-multigrid,
and HS-FD, the two-level hybrid Schwarz method with fast-diagonalization cell solves. Both build on
established methods. The Chebyshev smoother, the $p$-multigrid hierarchy, the hybrid Schwarz smoother with
FD local solves and the coarse corrections come from the literature. This work adds their
form on hanging-node octrees and distributed GPUs, the replicated and exact FD coarse solves, the
treatment of a variable coefficient and of stretched cells and inclusions, and the automatic choice
of parameters (\cref{tab:precond-parts}). The parts shared by both preconditioners are discussed first: the
Chebyshev smoother, the trilinear coarse level and the conditions under which a V-cycle is
symmetric positive definite (\cref{sec:precond-common}). The implementation on distributed GPUs,
which every solver of this work uses, and the precision of each part follow in
\cref{sec:implementation}. \Cref{sec:cjpmg,sec:asfd} then describe what is specific to
CJ-$p$MG and HS-FD, and \cref{sec:sm-precond} lists the solver settings.

\subsection{Chebyshev smoother, coarse level and symmetry}
\label{sec:precond-common}

\paragraph{Smoother}
Both preconditioners smooth with a Chebyshev polynomial in a preconditioned
operator~\cite{Adams2003PolynomialSmoothing}. For a symmetric positive
definite $\bB_{\mathrm{s}}$ and a window $[\lambda_0, \lambda_1]$, let $\theta =
(\lambda_1+\lambda_0)/2$, $\delta = (\lambda_1-\lambda_0)/2$ and $\tau_\nu(t) =
T_\nu\bigl((\theta-t)/\delta\bigr)/T_\nu(\theta/\delta)$, with $T_\nu$ the Chebyshev polynomial of
the first kind of degree $\nu$. The degree-$\nu$ smoother
\begin{equation}
  \mathcal{S}_\nu(\bB_{\mathrm{s}}) = q_\nu(\bB_{\mathrm{s}}^{-1}\bA)\,\bB_{\mathrm{s}}^{-1},
  \qquad q_\nu(t) = \frac{1-\tau_\nu(t)}{t},
  \label{eq:cheby}
\end{equation}
has the error propagator $\tau_\nu(\bB_{\mathrm{s}}^{-1}\bA)$ and, started from zero, costs
$\nu$ applications of $\bB_{\mathrm{s}}^{-1}$ and $\nu-1$ operator actions. The upper end of
every window is $\lambda_1 = 1.1\lambda_{\max}$, where $\lambda_{\max}$ is the largest Ritz value of a Lanczos
run on $\bB_{\mathrm{s}}^{-1}\bA$ at setup~\cite{Saad2003IterativeMethods}, started from a vector
that does not depend on the partition.

\paragraph{Coarsest level}
The coarsest level of both preconditioners is the space of trilinear elements on the cells of
the fine mesh, with one DoF per unconstrained vertex. The hanging-node, periodic and Dirichlet
constraints are imposed on this space as on the fine level. The dimension $n_c$ of the trilinear
space is a small fraction of $N$.

In the $p$-multigrid of Munch et al.\ and of nekRS the coarse level is solved by AMG over all MPI
ranks~\cite{Munch2023}. In this work the coarse system is instead replicated on every MPI rank, as in the redundant coarse solves common in earlier codes~\cite{Deville2002HighOrder}. One Allreduce per V-cycle sums the coarse residual onto every rank, and each rank
then solves the same coarse system, with no coarse ghost exchange and no distributed coarse
solve. Replication is chosen because the coarse level is small and because it avoids the setup of
an AMG hierarchy, which takes about three quarters of the time to solution of low-order-refined
AMG for a screened Poisson problem at $p = 6$ on one GPU~\cite{Pazner2023GPULOR}. The price is one reduction per V-cycle,
whose latency grows with the number of ranks, and a coarse solve repeated on every rank. A coarse
system solved by BoomerAMG is therefore distributed above $10^{6}$ unknowns
(\cref{sec:sm-auto}). The exact FD coarse solve can also be distributed without a global
reduction (\cref{sec:fd-coarse}), as on the 3D-IC rows on Aurora.

The coarse matrix is assembled on the trilinear space with the coefficient at the quadrature
points of the fine operator. The coarse system is solved in one of four ways, chosen by the rule
of \cref{sec:sm-auto}. On uniform and graded tensor meshes the trilinear grid is a tensor grid. Where the coefficient
varies at most along $z$, cell layer by cell layer, the coarse system is then solved exactly by
the FD solve of \cref{sec:fd-coarse}. Otherwise,
for $n_c \le 8000$, the explicit inverse is formed once on the GPU by a Cholesky factorization
and applied as a matrix--vector product. For larger $n_c$ the coarse system is solved by a
Jacobi-preconditioned Chebyshev iteration on the replicated sparse matrix or by BoomerAMG
V-cycles~\cite{Henson2002BoomerAMG} on the GPU, called from the hypre
library~\cite{Falgout2002hypre}. The Cholesky inverse,
the Chebyshev iteration and BoomerAMG are standard solvers. The replication, the exact FD coarse
solve and the automatic choice among the four are added in this work. Except for a distributed
BoomerAMG level, every coarse solver acts on the replicated system, so the choice adds no
communication. On the fully periodic operator the coarse solve acts on the mean-free subspace.

\paragraph{Symmetry}
CG needs a symmetric positive definite preconditioner. In every CG iteration the V-cycle is
applied to the residual $\br$ and returns the preconditioned residual $\bz = \bB^{-1}\br$. With
fixed smoother degrees, equal pre- and post-smoothing, a zero initial guess, restriction by the
transpose $\bP_\ell^{T}$ of the prolongation and a fixed linear coarse solve, $\bB^{-1}$ is a
linear map that does not change between iterations. Its error propagator is
$\bE = (\bI-\mathcal{S}\bA)(\bI-\bP\bB_c^{-1}\bP^{T}\bA)(\bI-\mathcal{S}\bA)$ and
$\bB^{-1} = (\bI-\bE)\bA^{-1}$, where $\bB_{\mathrm{s}}$ is the preconditioner inside the Chebyshev
smoother $\mathcal{S}$, the diagonal $\bD_\ell$ for CJ-$p$MG and the Schwarz operator for HS-FD, and
$\bB_c^{-1}$ is the approximate inverse on the coarse level. All three factors of $\bE$ are self-adjoint
in the $\bA$-inner product, so $\bB^{-1}$ is symmetric, as in the argument of Tatebe for multigrid
preconditioners with equal pre- and post-smoothing~\cite{Tatebe1993MGCG,Tatebe1994MGCGDistributed}.
It is positive definite when $\|\bE\|_{\bA} < 1$, which holds if (i) the eigenvalues of
$\bB_{\mathrm{s}}^{-1}\bA$ lie in $(0, \lambda_0 + \lambda_1)$, where the Chebyshev factor
$\tau_\nu$ of \cref{eq:cheby} has modulus below one, and (ii) the eigenvalues of
$\bB_c^{-1}\bP^{T}\bA\bP$ lie in $(0, 2]$, so that the coarse correction does not increase the error in the $\bA$-norm.

With $\lambda_1 = 1.1\lambda_{\max}$, condition (i) holds as long as the largest eigenvalue exceeds the
estimate $\lambda_{\max}$ by less than $0.1\lambda_{\max} + \lambda_0$. For
CJ-$p$MG the conditions apply level by level, with $\bB_c^{-1}$ the cycle of the next level,
and for HS-FD with $\bB_c^{-1}$ the coarse solve, $\bA_c^{-1}$ when it is exact (\cref{sec:asfd}). With a variable coefficient the
rediscretized level operators differ from the Galerkin products, and (ii) is assumed
(\cref{sec:sm-spd}).

\subsection{Implementation on distributed GPUs}
\label{sec:implementation}

The direct FD solve, CJ-$p$MG, HS-FD and Jacobi PCG share one
matrix-free implementation, written in CUDA for the A100 GPUs
of Polaris, in HIP for the MI250X of Frontier and in SYCL for the Max 1550 of Aurora,
with the same algorithms in each.

The fine operator is applied by sum factorization on the GLL nodes, in the fused kernels of
Kodali et al.~\cite{KodaliSC26}. Hanging-node constraints are resolved cell-locally inside every
kernel~\cite{Ljungkvist2017Matrix-FreeMeshes}. Ghost values are exchanged with GPU-aware MPI or the
vendor collective libraries on Polaris and Frontier and through host memory on Aurora, except for the
larger messages of the operator action on the 3D-IC stacks.
\Cref{sec:sm-impl} gives the kernels, exchanges and work per iteration.

\paragraph{Precision and parameter selection}
The outer CG, its vectors and all reductions are FP64 for every method, and the direct solve is
FP64 throughout. Inside the preconditioners the operator actions and the box solves use FP32
arithmetic, and they switch to FP64 when the largest cell aspect ratio or the largest ratio of
the coefficient within one cell is at least 300 (\cref{sec:sm-impl}). No parameter of either
preconditioner is set per test case. The Chebyshev windows, the coarse solver (exact FD, a dense
inverse, a Chebyshev iteration or BoomerAMG), the nodal scaling, the dense cut cells, the coarse
level of HS-FD and the precision follow from Lanczos estimates and short measured trials at setup
(\cref{sec:sm-auto}).

\subsection{Chebyshev--Jacobi \texorpdfstring{$p$}{p}-multigrid}
\label{sec:cjpmg}

CJ-$p$MG is the $p$-multigrid of R{\o}nquist and Patera~\cite{Ronquist1987SEMMultigrid} with
Chebyshev--Jacobi
smoothing~\cite{Sundar2015MultigridComparison}
and the smoothing band of nekRS and Ratel~\cite{Phillips2022TuningSEM}. Munch
et al.\ apply it on hanging-node octrees on CPUs~\cite{Munch2023}. This work runs every level,
with the hanging-node constraints of each level, on distributed GPUs, and adds the replicated
coarse level of \cref{sec:precond-common} and a lower end of the band for strongly stretched
cells.

\paragraph{Levels and cycle}
The levels are the continuous FE spaces of degrees $p_0 = p > p_1 > \dots > p_{L} = 1$ on the
fine mesh, with
\begin{equation}
  p_{\ell+1} = \lceil p_\ell/2 \rceil,
  \label{eq:pmg-degrees}
\end{equation}
so that $p = 7$ gives the degree sequence $7 \to 4 \to 2 \to 1$ on the Laplace and the
generalized Poisson operators. Going from the fine level straight to the trilinear level raises
the iterations at $p = 7$ threefold, from 6 to 18 on the uniform box of $20^3$ cells of \cref{tab:results-psweep}, so the
intermediate levels carry much of the error reduction.
Every level uses the cells, the partition and the ghost pattern of the fine level. The
hanging-node, periodic and Dirichlet constraints are imposed on each level separately. The
operator $\bA_\ell$ of level $\ell$ is obtained by discretizing the same problem again at degree
$p_\ell$ on the same cells, with the coefficient evaluated at the quadrature points of that
level. It is applied matrix-free by the fused sum-factorization kernels of the fine operator
(\cref{sec:implementation}). In every CG iteration one V(1,1) cycle approximately solves the
error equation for the current residual, starting from zero (\cref{alg:pmg-vcycle}). On each
level above the coarsest the cycle smooths, restricts the residual to the next level by
$\bP_\ell^{T}$, adds the correction $\bP_\ell \bz_{\ell+1}$ prolonged from the next level and
smooths again. The coarsest level is the replicated trilinear level of
\cref{sec:precond-common} (\cref{fig:pmg-vcycle}).

\paragraph{Smoother}
On each level above the coarsest the smoother is the Chebyshev polynomial \cref{eq:cheby} with
$\bB_{\mathrm{s}} = \bD_\ell$, the diagonal of $\bA_\ell$. The band is $[0.1\lambda_{\max},
1.1\lambda_{\max}]$, where $\lambda_{\max}$ is the largest Ritz value of $\bD_\ell^{-1}\bA_\ell$,
computed once per level, with a smaller lower end on cells of large aspect ratio
(\cref{sec:sm-auto}).
The polynomial degree is $\nu_0 = 6$ on the fine level for both operators. On the intermediate
levels it is $\nu_\ell = 3$ for the Laplace operator and $\nu_\ell = 4$ for the generalized
Poisson operator. On hanging-node meshes $\bD_\ell$ is the exact diagonal of the constrained
operator $\bC_\ell^{T}\bA_\ell\bC_\ell$, assembled cell by cell.

\paragraph{Transfers}
The prolongation $\bP_\ell$ interpolates a function of degree $p_{\ell+1}$ to degree $p_\ell$,
and the restriction is its transpose. On each cell $\bP_\ell$ is the tensor product of one
$(p_\ell+1)\times(p_{\ell+1}+1)$ interpolation matrix per direction and is applied by sum
factorization, as the operator is. Both transfers resolve the hanging-node constraints cell by
cell, as the operator kernels do, for every pair of levels, the step to the trilinear level
included. At setup each matrix-free transfer is compared with an assembled sparse transfer on a
test vector, and the sparse transfer is used instead if the two differ by more than a relative
$10^{-12}$.

\subsection{Two-level hybrid Schwarz with fast-diagonalization cell solves}
\label{sec:asfd}

HS-FD is the hybrid Schwarz--multigrid method of Lottes and
Fischer~\cite{Lottes2005HybridMultigrid}, an additive Schwarz smoother of FD local solves combined
multiplicatively with a coarse correction, with the Chebyshev acceleration of Phillips
et al.~\cite{Phillips2022TuningSEM}. Here the local problems are restricted to the nodes of one
cell, which keeps them Kronecker sums on 2:1 hanging-node octrees, the weights and local problems
are extended to a variable coefficient, and coarse corrections are added for strongly stretched
cells and isolated inclusions, all on distributed GPUs (\cref{tab:precond-parts}).

\paragraph{Cell boxes}
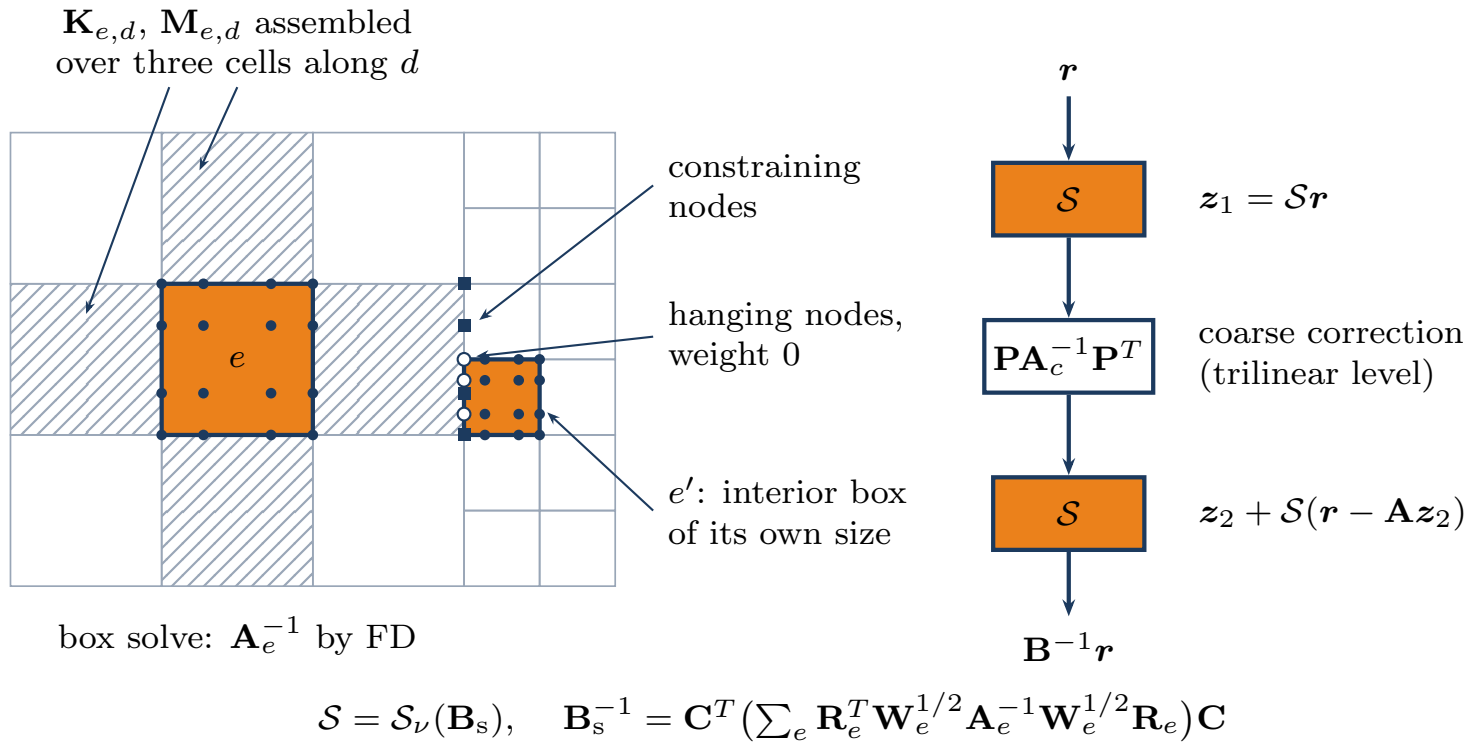
\begin{figure}[tbhp]
\centering
\scalebox{0.8}{\def\hsgll{-1,-0.4472135955,0.4472135955,1}%
\def\hsnodes#1#2#3{\foreach \u in \hsgll{\foreach \v in \hsgll{%
  \fill[deepnavy] ({#1+#3*(1+\u)/2},{#2+#3*(1+\v)/2}) circle (1pt);}}}%
\begin{tikzpicture}[font=\footnotesize,every node/.append style={text=black},
    hscell/.style={draw=deepnavy!45,line width=0.4pt},
    hsflow/.style={-{Stealth[length=4pt]},deepnavy,line width=1pt},
    hsnote/.style={-{Stealth[length=3pt]},deepnavy,line width=0.5pt},
    hsstep/.style={draw=deepnavy,line width=1pt,minimum width=1.25cm,minimum height=0.6cm,inner sep=2pt}]
  \begin{scope}[scale=1.25]
  \foreach \cx/\cy in {0/1,2/1,1/0,1/2}{%
    \fill[pattern=north east lines,pattern color=deepnavy!45] (\cx,\cy) rectangle ++(1,1);}
  \fill[mforange] (1,1) rectangle (2,2);
  \foreach \i in {0,1,2}{\foreach \j in {0,1,2}{\draw[hscell] (\i,\j) rectangle ++(1,1);}}
  \foreach \i in {0,1}{\foreach \j in {0,...,5}{\draw[hscell] ({3+0.5*\i},{0.5*\j}) rectangle ++(0.5,0.5);}}
  \fill[mforange] (3,1) rectangle (3.5,1.5);
  \draw[hscell] (3,1) rectangle (3.5,1.5);
  \draw[deepnavy,line width=1pt] (1,1) rectangle (2,2);
  \draw[deepnavy,line width=1pt] (3,1) rectangle (3.5,1.5);
  \hsnodes{1}{1}{1}
  \node[fill=mforange,inner sep=1pt] at (1.5,1.5) {$e$};
  \foreach \u in {-0.4472135955,0.4472135955,1}{\foreach \v in \hsgll{%
    \fill[deepnavy] ({3+0.25*(1+\u)},{1+0.25*(1+\v)}) circle (1pt);}}
  \fill[deepnavy] (3,1) circle (1pt);
  \foreach \v in {-0.4472135955,0.4472135955,1}{%
    \draw[deepnavy,fill=white,line width=0.5pt] (3,{1+0.25*(1+\v)}) circle (1.2pt);}
  \foreach \v in \hsgll{\fill[deepnavy] ($(3,{1+0.5*(1+\v)})+(-1.2pt,-1.2pt)$) rectangle ++(2.4pt,2.4pt);}
  \node[anchor=south,align=center] (asm) at (1.5,3.25) {$\bK_{e,d}$, $\bM_{e,d}$ assembled\\ over three cells along $d$};
  \draw[hsnote] ($(asm.south)+(-0.45,0.05)$) -- (0.5,1.8);
  \draw[hsnote] ($(asm.south)+(0.0,0.05)$) -- (1.25,2.75);
  \node[anchor=west,align=left] (cons) at (4.25,2.65) {constraining\\ nodes};
  \draw[hsnote] (cons.west) -- (3.1,1.76);
  \node[anchor=west,align=left] (hang) at (4.25,1.65) {hanging nodes,\\ weight $0$};
  \draw[hsnote] (hang.west) -- (3.08,1.52);
  \node[anchor=west,align=left] (fine) at (4.25,0.5) {$e'$: interior box\\ of its own size};
  \draw[hsnote] (fine.west) -- (3.55,1.15);
  \node[anchor=north,align=center] at (1.5,-0.1) {box solve: $\bA_e^{-1}$ by FD};
  \end{scope}
  \begin{scope}[shift={(8.75,0)}]
    \node[hsstep,fill=mforange] (s1) at (0,3.2) {$\mathcal{S}$};
    \node[hsstep,fill=white] (cc) at (0,1.9) {$\bP\bA_c^{-1}\bP^{T}$};
    \node[hsstep,fill=mforange] (s2) at (0,0.6) {$\mathcal{S}$};
    \draw[hsflow] (0,4.05) node[above] {$\br$} -- (s1);
    \draw[hsflow] (s1)--(cc); \draw[hsflow] (cc)--(s2);
    \draw[hsflow] (s2) -- (0,-0.25) node[below] {$\bB^{-1}\br$};
    \node[anchor=west] at (0.95,3.2) {$\bz_1 = \mathcal{S}\br$};
    \node[anchor=west,align=left] at (0.95,1.9) {coarse correction\\ (trilinear level)};
    \node[anchor=west] at (0.95,0.6) {$\bz_2 + \mathcal{S}(\br - \bA\bz_2)$};
  \end{scope}
  \node[anchor=north] at (6.3,-0.75)
    {$\mathcal{S} = \mathcal{S}_\nu(\bB_{\mathrm{s}})$, \quad
     $\bB_{\mathrm{s}}^{-1} = \bC^{T}\bigl(\sum_e \bR_e^{T}\bW_e^{1/2}\bA_e^{-1}\bW_e^{1/2}\bR_e\bigr)\bC$};
\end{tikzpicture}}
\caption[Cell boxes and cycle of HS-FD]{HS-FD on a two-dimensional slice at FE degree 3. The box
of cell $e$ holds the nodes of $e$ and couples to its neighbors through one-dimensional matrices
assembled over three cells along each direction. At a 2:1 interface the cell $e'$ receives the
interior box of a cell of its own size, with zero weight on its hanging nodes. The weighted box
solves form the smoother $\mathcal{S}$ of the two-level cycle with the trilinear coarse
correction.}
\label{fig:asfd-boxes}
\end{figure}
Lottes and Fischer restrict $\bA$ to one element and one further layer of nodes inside each
neighbor~\cite{Lottes2005HybridMultigrid}. HS-FD restricts $\bA$ to
the nodes of one cell, so neighboring boxes share only their common face (\cref{fig:asfd-boxes}). For a cell $e$ and a direction $d \in \{x,y,z\}$, let $\bK_{e,d}$ and
$\bM_{e,d}$ be the $(p+1)\times(p+1)$ blocks on the nodes of $e$ of the one-dimensional stiffness
and mass matrices assembled over three cells of the size of $e$, with $e$ in the middle. They use
the Gauss quadrature of the operator, so the mass matrices are consistent. The box is their
Kronecker sum,
\begin{equation}
  \bA_e = \bK_{e,x}\otimes\bM_{e,y}\otimes\bM_{e,z} + \bM_{e,x}\otimes\bK_{e,y}\otimes\bM_{e,z}
        + \bM_{e,x}\otimes\bM_{e,y}\otimes(\bK_{e,z} + \sigma\bM_{e,z}).
  \label{eq:asfd-box}
\end{equation}
Where the neighbors of $e$ have its size, $\bA_e = \bR_e\bA\bR_e^{T}$ is the block of the
assembled operator on the $(p+1)^3$ nodes of the cell, with $\bR_e$ the restriction to those
nodes. Unlike the Neumann matrix of the cell alone, the box includes the coupling to the neighbors and
is symmetric positive definite. At a Dirichlet face the boundary node is decoupled, and across a
periodic face the box of an interior cell is used. As in
\cref{sec:fd-inverse}, the inverse of the box follows from the generalized eigenpairs of its
one-dimensional matrices and is applied exactly in $\order(p^4)$ operations,
\begin{equation}
  \begin{aligned}
  &\bK_{e,d}\bS_{e,d} = \bM_{e,d}\bS_{e,d}\bLam_{e,d}, \qquad
  \bS_{e,d}^{T}\bM_{e,d}\bS_{e,d} = \bI,\\
  &\bA_e^{-1} = \bS_e\,
    \diag\Bigl(\frac{1}{\lambda_{e,x,i}+\lambda_{e,y,j}+\lambda_{e,z,k}+\sigma}\Bigr)\,\bS_e^{T},
  \end{aligned}
  \label{eq:asfd-apply}
\end{equation}
where $\bS_e = \bS_{e,x}\otimes\bS_{e,y}\otimes\bS_{e,z}$. For a brick cell of edge lengths
$h_d$ the eigenpairs follow from three reference boxes per direction (lower face, interior, upper
face), with eigenvalues scaled by $h_d^{-2}$ and eigenvectors by $h_d^{-1/2}$, so graded tensor
meshes need no further eigenproblems. On the GPU the box solve is a sum-factorized kernel with
$\bS_{e,d}$ in place of the shape-function matrices (\cref{sec:sm-impl}).

\paragraph{Weighted smoother}
The cell inverses are summed with two-sided weights,
\begin{equation}
  \bB_{\mathrm{s}}^{-1} = \bC^{T}
    \Bigl(\sum_e \bR_e^{T}\bW_e^{1/2}\bA_e^{-1}\bW_e^{1/2}\bR_e\Bigr)\bC,
  \label{eq:asfd-smoother}
\end{equation}
where $\bC$ maps the free DoFs to all mesh nodes (\cref{sec:problem}) and the diagonal matrix
$\bW_e$ holds the weights of the nodes of cell $e$. At every unconstrained node the weights of the
cells sharing it sum to one. The two-sided form of Lottes and Fischer keeps $\bB_{\mathrm{s}}^{-1}$
symmetric, so the preconditioner can be used inside CG (\cref{sec:sm-precond}). Their counting
weight, the inverse of the number of boxes that contain a node, counted over all MPI ranks, is
used here for a constant coefficient. With a variable coefficient the largest eigenvalue of
$\bB_{\mathrm{s}}^{-1}\bA$ grows in proportion to the contrast between neighboring cells,
as about 0.27 times the contrast between $10^2$ and $10^4$, and the Chebyshev window leaves a growing
part of the spectrum to the coarse level. For the generalized Poisson operator the weight of a
node in cell $e$ is its diagonal entry in the element matrix of $e$ divided by its
diagonal entry in $\bA$, normalized to one at cell-interior nodes. At a constrained node
the sum of its element entries replaces its entry in $\bA$.
\paragraph{Hanging nodes}
On an octree with 2:1 hanging nodes, a cell at a refinement interface receives the box of an
interior cell of its own size. Its hanging nodes receive zero weight and drop out of the
smoother, while its other nodes, the constraining nodes included, keep positive weights. The
smoother stays symmetric positive definite: its quadratic form is a sum of nonnegative cell
terms, and every free node has a positive weight in some cell. The approximation at the
interface can therefore cost iterations but not the convergence of CG. Resolving the hanging
nodes inside the box, as the operator kernel does, also gives a symmetric smoother, but with a larger maximum eigenvalue of $\bB_{\mathrm{s}}^{-1}\bA$ and more iterations.

\paragraph{Chebyshev acceleration}
As proposed by Phillips et al.~\cite{Phillips2022TuningSEM}, the Schwarz smoother is accelerated
by the Chebyshev polynomial \cref{eq:cheby} with $\bB_{\mathrm{s}}$ of \cref{eq:asfd-smoother}.
Adams et al.\ and Phillips et al.\ set the lower end $\lambda_0$ of the window to a fixed fraction
of $\lambda_{\max}$~\cite{Adams2003PolynomialSmoothing,Phillips2022TuningSEM}, and Phillips and Fischer
derive it from a bound on the V-cycle~\cite{Phillips2025OptimalSmoothers}. Here the selection of
\cref{sec:sm-auto} chooses between a fixed fraction of $\lambda_{\max}$ and the smallest Ritz value of a
Lanczos run on the smoothed operator restricted to the $\bA$-orthogonal complement of the coarse
space, the part of the spectrum the smoother must reduce, which is a candidate only when the coarse
solve is accurate enough for that run. The degree $\nu$ is 5 for $p < 6$. At
$p \ge 6$ it is 4 for the Laplace operator and, for the generalized Poisson operator, 5 with the
nodal scaling described below or 6 without it (\cref{sec:sm-precond}).

\paragraph{Two-level cycle}
With $\mathcal{S} = \mathcal{S}_\nu(\bB_{\mathrm{s}})$, the preconditioner is one multiplicative
V(1,1) cycle on $\bA\bz = \br$ from a zero initial guess,
\begin{equation}
  \bz_1 = \mathcal{S}\br, \qquad
  \bz_2 = \bz_1 + \bP\bA_c^{-1}\bP^{T}(\br - \bA\bz_1), \qquad
  \bB^{-1}\br = \bz_2 + \mathcal{S}(\br - \bA\bz_2),
  \label{eq:asfd-cycle}
\end{equation}
where $\bA_c$ is the matrix of the trilinear coarse level of \cref{sec:precond-common} and $\bP$
the interpolation from the trilinear space. The interpolation is applied matrix-free cell by
cell as a tensor product of one-dimensional matrices. Unlike CJ-$p$MG, HS-FD needs
no intermediate levels, because its cell solves smooth far more strongly than the point-Jacobi
Chebyshev iteration.
By conditions (i) and (ii) of \cref{sec:precond-common}, with $\bB_c^{-1} = \bA_c^{-1}$ for an exact coarse solve, $\bB^{-1}$ is symmetric positive definite when
the eigenvalues of $\bB_{\mathrm{s}}^{-1}\bA$ lie below $\lambda_0+\lambda_1$ and those of
$\bA_c^{-1}\bP^{T}\bA\bP$ in $(0,2]$, so the outer method is CG. A coarse system that is not
solved exactly is solved by the rule of \cref{sec:sm-auto}.

\paragraph{Coarse correction on stretched cells}
On strongly stretched cells, the trilinear coarse level does not represent the error that varies
slowly along the long direction of the cells, and HS-FD converges slowly
(\cref{sec:results-apps}). Where the largest cell aspect ratio is 4 or more, the setup also builds one BoomerAMG V-cycle on
the low-order-refined operator of the fine mesh~\cite{Pazner2020LOR}, assembled on the GPU, as a
second candidate for $\bP\bA_c^{-1}\bP^{T}$ in \cref{eq:asfd-cycle}, and the trial of
\cref{sec:sm-auto} keeps the faster one. Condition (ii) of \cref{sec:precond-common} is assumed for it, as for
the rediscretized levels.

\paragraph{Coarse vectors for isolated inclusions}
Each isolated inclusion of high coefficient, as in the battery electrode
(\cref{tab:results-coefficients}), adds an error component of small energy that is nearly
constant on the inclusion, and the trilinear coarse level does not represent it either. Tang and
Vuik~\cite{Tang2007LevelsetDeflation} remove such components in CG on bubbly flow by deflation with
one vector per bubble that equals one on it. The GenEO coarse spaces of Spillane et
al.~\cite{Spillane2014GenEO} treat general coefficient jumps in Schwarz methods through local
generalized eigenproblems. HS-FD adds vectors of the first kind automatically: when the coefficient ratio
exceeds $10^{2}$ and the trial of \cref{sec:sm-auto} converges slowly, it adds one coarse vector for
each inclusion that touches no Dirichlet face, equal to one on the inclusion and extended
harmonically into one surrounding layer of cells. These vectors enter the coarse correction additively,
through their Galerkin matrix with $\bA$, and CJ-$p$MG does not use them (\cref{sec:sm-precond}).

\paragraph{Generalized Poisson operator}
For the generalized Poisson operator $\bK^{\epsilon}$ the variable coefficient $\epsilon$ breaks
the Kronecker structure of the box. In this work $\epsilon$ is fitted on each cell by a separable
function $a(x)\,b(y)\,c(z)$ with affine logarithm. The weighted stiffness then factors into
one-dimensional matrices weighted by $a$, $b$ and $c$, so the fitted box is again a Kronecker sum
and \cref{eq:asfd-apply} inverts it exactly. For discontinuous Galerkin methods, Pazner and
Persson~\cite{Pazner2018TensorPrecond} and Diosady and Murman~\cite{Diosady2019TensorPrecond} fit
the element matrix algebraically by Kronecker products of one-dimensional matrices. Here the
coefficient is fitted instead, and the box is assembled from the fit.
Where the fit deviates from the coefficient, a nodal scaling enters on both sides of the boxes
through $\bW^{1/2}$, which keeps the smoother symmetric. On cells cut by an interface between two
materials, where the coefficient jumps, the natural logarithm of the coefficient deviates from the fit by
more than 1. There the box solve is replaced by an exact solve with the cell block
$\bR_e\bA\bR_e^{T}$ of $\bA$, factorized by Cholesky in FP64 and stored as its
weighted inverse. Dense cells are admitted within limits on cost and storage, and the selection of
\cref{sec:sm-auto} decides whether they are used, together with the nodal scaling
(\cref{sec:sm-precond}).
The
dense inverses are applied by matrix--vector products, to all dense cells in one kernel, in FP64.

\section{Numerical results}
\label{sec:results}

The results follow the cost model of \cref{sec:problem-cost}: the iterations
(\cref{sec:results-iterations}), the cost of one iteration (\cref{sec:results-cost}) and
their product, the time to solution (\cref{sec:results-time}). The later subsections treat blocks
of right-hand sides (\cref{sec:results-mrhs}), the application problems (\cref{sec:results-apps}),
the comparison with LOR-AMG and nekRS (\cref{sec:results-amg}) and the selection rule
(\cref{sec:results-selection}).

\subsection{Systems, machines and measurement basis}
\label{sec:results-setup}

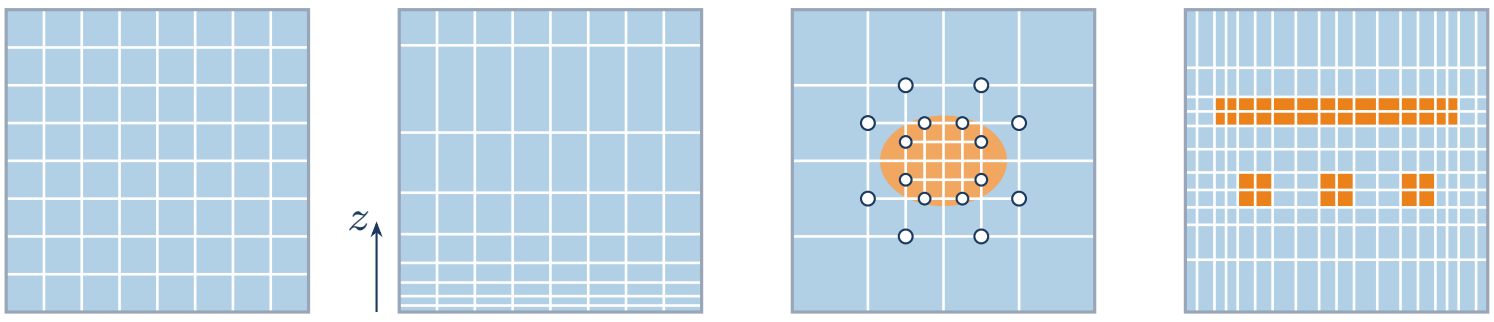
\begin{figure}[tbhp]
\centering
\def\mfmeshw{2.0}%
\begin{tikzpicture}[font=\footnotesize,
    mfmesh/.style={draw=white,line width=0.5pt},
    mfframe/.style={draw=deepnavy!45,line width=0.6pt},
    mfhang/.style={draw=deepnavy,fill=white,line width=0.4pt}]
  \begin{scope}
    \fill[mfblue!35] (0,0) rectangle (\mfmeshw,\mfmeshw);
    \foreach \t in {1,...,7}{\draw[mfmesh] ({\t*\mfmeshw/8},0) -- ++(0,\mfmeshw); \draw[mfmesh] (0,{\t*\mfmeshw/8}) -- ++(\mfmeshw,0);}
    \draw[mfframe] (0,0) rectangle (\mfmeshw,\mfmeshw);
    \node[below] at ({\mfmeshw/2},0) {(a) uniform};
  \end{scope}
  \begin{scope}[xshift=2.6cm]
    \fill[mfblue!35] (0,0) rectangle (\mfmeshw,\mfmeshw);
    \foreach \t in {1,...,7}{\draw[mfmesh] ({\t*\mfmeshw/8},0) -- ++(0,\mfmeshw);}
    \foreach \z in {0.043,0.105,0.195,0.325,0.515,0.789,1.187,1.764}{\draw[mfmesh] (0,\z) -- ++(\mfmeshw,0);}
    \draw[mfframe] (0,0) rectangle (\mfmeshw,\mfmeshw);
    \draw[-{Stealth[length=3pt]},deepnavy] (-0.15,0) -- (-0.15,0.6) node[left,inner sep=1pt] {$z$};
    \node[below] at ({\mfmeshw/2},0) {(b) graded};
  \end{scope}
  \begin{scope}[xshift=5.2cm]
    \fill[mfblue!35] (0,0) rectangle (\mfmeshw,\mfmeshw);
    \fill[mforange!70] ({\mfmeshw/2},{\mfmeshw/2}) ellipse (0.42 and 0.3);
    \foreach \t in {1,2,3}{\draw[mfmesh] ({\t*\mfmeshw/4},0) -- ++(0,\mfmeshw); \draw[mfmesh] (0,{\t*\mfmeshw/4}) -- ++(\mfmeshw,0);}
    \foreach \t in {3,5}{\draw[mfmesh] ({\t*\mfmeshw/8},{\mfmeshw/4}) -- ++(0,{\mfmeshw/2}); \draw[mfmesh] ({\mfmeshw/4},{\t*\mfmeshw/8}) -- ++({\mfmeshw/2},0);}
    \foreach \t in {7,9}{\draw[mfmesh] ({\t*\mfmeshw/16},{3*\mfmeshw/8}) -- ++(0,{\mfmeshw/4}); \draw[mfmesh] ({3*\mfmeshw/8},{\t*\mfmeshw/16}) -- ++({\mfmeshw/4},0);}
    \draw[mfframe] (0,0) rectangle (\mfmeshw,\mfmeshw);
    \foreach \t in {3,5}{%
      \draw[mfhang] ({\mfmeshw/4},{\t*\mfmeshw/8}) circle (1.3pt); \draw[mfhang] ({3*\mfmeshw/4},{\t*\mfmeshw/8}) circle (1.3pt);
      \draw[mfhang] ({\t*\mfmeshw/8},{\mfmeshw/4}) circle (1.3pt); \draw[mfhang] ({\t*\mfmeshw/8},{3*\mfmeshw/4}) circle (1.3pt);}
    \foreach \t in {7,9}{%
      \draw[mfhang] ({3*\mfmeshw/8},{\t*\mfmeshw/16}) circle (1.1pt); \draw[mfhang] ({5*\mfmeshw/8},{\t*\mfmeshw/16}) circle (1.1pt);
      \draw[mfhang] ({\t*\mfmeshw/16},{3*\mfmeshw/8}) circle (1.1pt); \draw[mfhang] ({\t*\mfmeshw/16},{5*\mfmeshw/8}) circle (1.1pt);}
    \node[below] at ({\mfmeshw/2},0) {(c) octree};
  \end{scope}
  \begin{scope}[xshift=7.8cm,scale=0.7692]
    \fill[mfblue!35] (0,0) rectangle (2.6,2.6);
    \foreach \x in {0.45,1.15,1.85}{\fill[mforange] (\x,0.9) rectangle ++(0.3,0.3);}
    \fill[mforange] (0.25,1.6) rectangle (2.35,1.85);
    \foreach \x in {0.100,0.250,0.350,0.450,0.600,0.750,0.950,1.150,1.300,1.450,1.650,1.850,2.000,2.150,2.250,2.350,2.500}{\draw[mfmesh] (\x,0) -- ++(0,2.6);}
    \foreach \z in {0.450,0.750,0.900,1.050,1.200,1.400,1.600,1.725,1.850,2.100}{\draw[mfmesh] (0,\z) -- ++(2.6,0);}
    \draw[mfframe] (0,0) rectangle (2.6,2.6);
    \node[below] at ({2.6/2},0) {(d) interconnect};
  \end{scope}
\end{tikzpicture}
\caption[Meshes of the test problems]{Schematic two-dimensional slices of the meshes, not renders of the
meshes used: (a) uniform (U-D, U-PD, U-P, uniform box), (b) graded along $z$ (G, G($R$)), (c)
octree with 2:1 hanging nodes (circles) around an inner region (H-GP, H-GP-s), (d) conductors in a
graded tensor-product mesh (interconnects of \cref{tab:results-apps}).}
\label{fig:results-meshes}
\end{figure}

\Cref{tab:results-iterations} lists the test problems, and \cref{tab:sm-difficulty} gives for each
one the property that makes it hard for a solver, such as strong grading or a jump in the
coefficient. The problems lie on the meshes sketched in \cref{fig:results-meshes}. Each problem is named by its
mesh and boundary conditions. U-D, U-PD and U-P are uniform meshes with Dirichlet conditions, with
periodic conditions in $x$ and $y$, and with periodic conditions in all three directions, and G is
a mesh graded along $z$. On these four the operator is the Laplacian $\bK$. H-GP and H-GP-s are
octree meshes with the generalized Poisson operator $\bK^{\epsilon}$, whose coefficient
$\epsilon(\bx)$ equals 1 inside an inner region and rises to $\epsilon_\infty = 78.4$
outside it across a thin layer, so that $\epsilon$ varies strongly inside the cells of that layer.
The FE degree is $p = 7$ unless stated otherwise, with the Gauss quadrature of
\cref{sec:problem-operators}. The meshes and the coefficient come from DFT-FE electrostatics of an aluminium
cluster, a garnet electrolyte, bulk and slab molybdenum, and a protein and a molecule in implicit
solvent (\cref{sec:sm-difficulty}).

\paragraph{Machines} Polaris (Argonne Leadership Computing Facility, ALCF) has four NVIDIA A100
GPUs (40\,GB) per node. Frontier (Oak Ridge Leadership Computing Facility) has four AMD MI250X
per node, each with two GCDs, and Aurora (ALCF) six Intel Data Center
GPU Max 1550 (PVC) per node, each with two tiles. One MPI rank drives one A100, one GCD or one
PVC tile. The compilers and libraries are listed in \cref{sec:sm-impl}.
Aurora is allocated in whole nodes of 12 tiles, so it uses whole-node counts near the rank counts
of Polaris and Frontier, such as 12 tiles for 8 ranks and 24 or 36 for 32.
Every table gives the ranks.

\paragraph{Measurement} Each configuration is solved six times in one run, each time from a zero
initial guess. The first solve is discarded and the tables report solves 2 to 6. Because atomic
updates and reductions on the GPU accumulate in an order that varies from one solve to the next,
the number of iterations can differ among these five solves. The tables give the number that
occurs most often, or the median if none does, with the range where the solves differ. The times
are medians over the same five solves, with the setup excluded and the internal phase timers
switched off. Conjugate gradient (CG) solver time is the wall time for the CG iterations in a single solve. Fast diagonalization (FD) time is the wall time of one distributed FD solve, including the redistribution into and out of the FD layout. If the
true residual of the first FD solve misses the tolerance, every FD solve adds one step of
iterative refinement in FP64, and the FD time includes it. For a block of $m$ right-hand sides
both times are the time of one block solve divided by $m$. Neither time includes the assembly of
the right-hand side and its copy to the device, which precede the solve, or the application of
the constraints to the converged solution and its copy to the host, which follow it.
The CG solver time also excludes the initial residual and the first preconditioner
application, which precede the iterations. The Chebyshev--Jacobi preconditioner
has degree 5 on $[\lambda_{\min}, \lambda_{\max}]$ from a 20-step Lanczos run. The CJ-$p$MG
smoothers have degree 6 on the fine level and 3 (Laplacian) or 4 (generalized Poisson
operator) on the intermediate levels. The Chebyshev windows, the eigenvalue estimates and the
coarse solver follow from the operator, with no setting chosen per case
(\cref{sec:implementation}).

\paragraph{Stopping test} Every solve stops at a relative true residual of $10^{-10}$, at
$10^{-9}$ on two problems that do not reach $10^{-10}$ in FP64, or after $20\,000$ iterations
(\cref{sec:sm-stopping}). At coefficient ratios of $10^{3}$ or more the results give the normwise backward error $\eta = \|\bff - \bA\bu\|_2 / (\|\bA\|_2 \|\bu\|_2 +
\|\bff\|_2)$, with $\|\bA\|_2$ from 20 Lanczos steps.

\paragraph{Comparison across machines} Iterations are compared across machines, and wall times
only within one machine and one rank count. Where the iterations on two machines differ by more
than one, the text names the cause where it is known. Such a difference can come from hypre, whose BoomerAMG
hierarchy differs between GPU architectures, as for HS-FD on the graded mesh G (18 iterations on
Aurora, 23 on Polaris and Frontier), or from the automatic selection, timed on each machine.

\paragraph{Comparator codes} The deal.II $p$-multigrid runs on the host CPU and is compared in
iterations only. LOR-AMG through MFEM is not run on Aurora, and nekRS is compared on all
three machines (\cref{sec:sm-comparators}).

\subsection{Iterations}
\label{sec:results-iterations}

\begin{table}[tbhp]
\caption[Test systems, time and iterations]{Test systems and time to solution, with the
iterations in parentheses, from a zero initial guess to a relative true residual of $10^{-10}$, on
the basis of \cref{sec:results-setup}, on Polaris (P, A100 GPUs), Frontier (F, MI250X GCDs) and
Aurora (A, Intel Max 1550 tiles, host-staged MPI): FD time of one solve in ms with its relative
true residual below it, CG solver time of one solve in s. FE degree 7. Systems, with $N$ and the
boundary conditions in $x$, $y$ and $z$ (D: Dirichlet, P: periodic): uniform tensor-product meshes
U-D 11\,390\,625 (DDD), U-PD 12\,525\,248 (PPD) and U-P 11\,390\,625 (PPP), the
tensor-product mesh G graded along $z$ with ratio $R = 2^{14}$ of the largest to the smallest cell
size, 11\,224\,473 (PPD), and octree meshes with hanging nodes H-GP 11\,390\,280 and H-GP-s
5\,417\,794 (DDD). Operators: Laplacian $\bK$, singular on U-P, and generalized Poisson
$\bK^{\epsilon}$ with $\epsilon_\infty = 78.4$ on H-GP, whose times are in
\cref{fig:results-trp-scaling}, and H-GP-s. Outer CG and FD solve in FP64,
preconditioners with FP32 arithmetic on the fine and intermediate levels and in the HS-FD box
solves and FP64 elsewhere (\cref{sec:implementation}), with these parts in FP64 on G.
Jacobi PCG uses two-reduction CG on U-P and G, all other entries single-reduction
CG, flexible for CJ-$p$MG and HS-FD. A dash marks a method that does not apply (FD needs a tensor-product mesh) or is not
reported. $^{a}$Stops at the 20\,000-iteration cap with a relative residual of the CG recurrence
of $2.4\times10^{-5}$. $^{b}$On Frontier the BoomerAMG setup of the low-order-refined coarse level
of HS-FD gives 18 or 23 to 24 iterations from one run to the next.}
\label{tab:results-iterations}
\label{tab:results-systems}
\label{tab:results-time}
\centering
\setlength{\tabcolsep}{3pt}
\begin{tabular}{@{}llrrrrr@{}}
\toprule
system & & ranks & FD & Jacobi & CJ-$p$MG & HS-FD \\
\midrule
U-D & P & 32 & \makecell[r]{0.785\\ $3.4\times10^{-12}$} & \makecell[r]{0.335\\ (1122)} & \makecell[r]{0.0253\\ (6)} & \makecell[r]{0.0192\\ (7)} \\
 & F & 32 & \makecell[r]{1.40\\ $2.8\times10^{-12}$} & \makecell[r]{0.253\\ (1122)} & \makecell[r]{0.0225\\ (6)} & \makecell[r]{0.0170\\ (7)} \\
 & A & 24 & \makecell[r]{1.54\\ $2.7\times10^{-12}$} & \makecell[r]{0.557\\ (1122)} & \makecell[r]{0.0731\\ (6)} & \makecell[r]{0.0703\\ (7)} \\
U-PD & P & 32 & \makecell[r]{0.766\\ $4.9\times10^{-13}$} & \makecell[r]{0.582\\ (2497)} & \makecell[r]{0.0247\\ (6)} & \makecell[r]{0.0177\\ (8)} \\
 & F & 32 & \makecell[r]{1.33\\ $4.7\times10^{-13}$} & \makecell[r]{0.657\\ (2497)} & \makecell[r]{0.0222\\ (6)} & \makecell[r]{0.0178\\ (8)} \\
 & A & 36 & \makecell[r]{1.24\\ $3.9\times10^{-13}$} & \makecell[r]{0.824\\ (2497)} & \makecell[r]{0.0477\\ (6)} & \makecell[r]{0.0464\\ (8)} \\
U-P & P & 32 & \makecell[r]{1.53\\ $1.5\times10^{-13}$} & \makecell[r]{0.759\\ (1628)} & \makecell[r]{0.0439\\ (7)} & \makecell[r]{0.0255\\ (7)} \\
 & F & 32 & \makecell[r]{1.67\\ $1.5\times10^{-13}$} & \makecell[r]{0.544\\ (1628)} & \makecell[r]{0.0306\\ (7)} & \makecell[r]{0.0188\\ (7)} \\
 & A & 36 & \makecell[r]{5.24\\ $1.3\times10^{-13}$} & \makecell[r]{0.781\\ (1628)} & \makecell[r]{0.0708\\ (7)} & \makecell[r]{0.0527\\ (7)} \\
G & P & 32 & \makecell[r]{1.10\\ $9.8\times10^{-13}$} & \makecell[r]{136.9\\ (349\,666--\\ 351\,501)} & \makecell[r]{0.970\\ (250)} & \makecell[r]{0.134\\ (23)} \\
 & F & 32 & \makecell[r]{1.40\\ $9.3\times10^{-13}$} & \makecell[r]{109.2\\ (349\,712--\\ 352\,326)} & \makecell[r]{0.945\\ (250)} & \makecell[r]{0.104$^{b}$\\ (23)} \\
 & A & 36 & \makecell[r]{1.31\\ $7.4\times10^{-13}$} & (${>}20\,000$)$^{a}$ & \makecell[r]{2.38\\ (250)} & \makecell[r]{0.334\\ (18)} \\
H-GP & P & 32 & --- & (2317--2320) & (9) & (11) \\
 & F & 32 & --- & (2317--2320) & (9) & (11) \\
 & A & 36 & --- & (2317--2319) & (10) & (11) \\
H-GP-s & P & 8 & --- & (1634--1635) & (8) & (12) \\
 & F & 8 & --- & (1634--1635) & (8) & (12) \\
 & A & 12 & --- & (1635) & (8) & (11) \\
\bottomrule
\end{tabular}
\end{table}

On Polaris and Frontier the iterations agree to within the spread of the repeated solves
(\cref{tab:results-iterations}). On Aurora, they are the same on the four Poisson problems U-D,
U-PD, U-P and G, except for HS-FD on G, which needs 18 iterations instead of 23, and Jacobi PCG on G, which stops at the 20\,000-iteration cap (note a). The direct FD solve reaches a relative true
residual about two to three orders of magnitude below the tolerance of the iterative methods.

On the Poisson problems, the iterations of Jacobi PCG grow with the condition number of $\bK$, up
to about 350\,000 on the graded mesh G ($R = 2^{14}$). There CJ-$p$MG needs 250 iterations, and HS-FD, with the coarse level of BoomerAMG on the low-order-refined operator chosen
automatically (\cref{sec:asfd}), needs 23 on Polaris and Frontier,
against 36 for LOR-AMG (\cref{sec:results-amg}). On Frontier HS-FD is 1.44 times faster than LOR-AMG and 9.1 times faster than CJ-$p$MG, and LOR-AMG
is 6.3 times faster than CJ-$p$MG.
For the generalized Poisson operator the cell problem of HS-FD replaces $\epsilon$ on each cell by
a product of one-dimensional functions (\cref{sec:asfd}), while on H-GP and H-GP-s $\epsilon$
varies strongly inside the cells of its transition layer. CJ-$p$MG uses on every level the exact
diagonal of $\bK^{\epsilon}$ with the hanging-node constraints applied
(\cref{sec:problem-constraints}), which follows the coefficient node by node.
On H-GP HS-FD therefore needs 11
iterations against 9 for CJ-$p$MG (10 on 36 Aurora tiles), and its cheaper cycle still gives the
shorter CG solver time on every machine
(\cref{sec:results-time}).
Chebyshev--Jacobi PCG uses a Chebyshev polynomial of the Jacobi-preconditioned operator alone as
the preconditioner of
CG~\cite{KodaliSC26}. On H-GP on 32 Frontier GCDs it needs 628 iterations and
0.72\,s per solve, against 2317--2320 iterations for Jacobi PCG and 9 for CJ-$p$MG.
Under uniform refinement at $p = 4$ the iterations of both preconditioners are independent of
$h$, 5 for CJ-$p$MG and 5 to 6 for HS-FD, while those of Jacobi PCG double as $h$ is halved.
Between $p = 2$ and 7 CJ-$p$MG needs 5 to 8 iterations, while those of HS-FD grow from 4 to 7 on a
uniform and from 6 to 11 on an octree mesh, as expected for cell subdomains with minimal
overlap (\cref{sec:asfd}). With the degree sequence and settings of deal.II, CJ-$p$MG needs 7,
7 and 8 iterations at $p = 2$, 4 and 7, against 6, 6 and 8 for its host $p$-multigrid
(\cref{sec:sm-psweep}).

On the grading series G($R$) from $R = 1$ to $15\,625$ the time of one FD solve on four A100 GPUs
stays between 0.248 and 0.259\,ms up to $R = 4096$ and takes 0.913\,ms at $R = 15\,625$, including
one step of iterative refinement, while Jacobi PCG grows from 605 to 54\,944 iterations. HS-FD
needs fewer iterations than CJ-$p$MG at every grading ratio from $R = 64$, 18 against 446 at
$R = 15\,625$ (\cref{sec:sm-grading}).

The exact FD solve also serves as the coarse solve of both preconditioners whenever the
$\mathcal{Q}_1$ coarse grid is a tensor-product grid. On U-PD it lowers the iterations from 12 to 6
for CJ-$p$MG and from 13 to 8 for HS-FD against the coarse solve that the automatic rule selects
otherwise (\cref{sec:sm-fdcoarse}).

\Cref{tab:results-coefficients} adds three problems with other coefficients: a droplet of one fluid
in another, where the coefficient is the inverse density and the density ratio is $10^{3}$, the
screened operator $\bK + \sigma\bM$, and a layered medium whose permeability spans a coefficient
ratio of $10^{6}$. On these problems CJ-$p$MG needs 4 to 10 iterations and HS-FD 2 to 10, against 31 to 44 for
LOR-AMG and 28 to 2384 for Jacobi PCG. On the layered medium Jacobi PCG converges only with
two-reduction CG.
In the layered medium the permeability $\kappa(z)$ is
piecewise constant with $\log_{10} \kappa$ uniform in $[-3, 3]$, so the operator remains a
Kronecker sum and the direct FD solve applies. On 8 Frontier GCDs it takes 0.527\,ms at $p = 4$ and
2.35\,ms at $p = 7$ (\cref{tab:results-amg}), the latter with one step of iterative refinement.

The battery electrode, the last problem of \cref{tab:results-coefficients}, is a generalized Poisson
problem (\cref{eq:prob-strong}) on a fully periodic cube. Its coefficient models a lithium-ion
electrode as 390 disjoint conducting spheres in a matrix whose conductivity is $10^{-3}$ or
$10^{-4}$ times that of the spheres, a coefficient ratio of $10^{3}$ or $10^{4}$. At the ratio
$10^{4}$ the true residual stops falling near $1.4\times10^{-10}$, so these solves stop at a
relative residual of $10^{-9}$. CJ-$p$MG needs 72
iterations at $10^{3}$ and 180 at $10^{4}$. With the coarse vectors of the isolated components (\cref{sec:asfd}) HS-FD needs 63 and
397 to 399. HS-FD is faster at $10^{3}$ and CJ-$p$MG at $10^{4}$ on all three machines
(\cref{sec:sm-battery}).

\subsection{Cost per iteration}
\label{sec:results-cost}

\label{sec:problem-cost}
The time to solution of a preconditioned Krylov method is written as the
number of iterations times the cost of one iteration,
\begin{equation}
  T_{\mathrm{solve}} = n_{\mathrm{it}}
  \left( a\, t_{\mathrm{op}} + t_{\mathrm{pre}} + t_{\mathrm{vec}}
         + t_{\mathrm{comm}} \right),
  \label{eq:prob-cost}
\end{equation}
where $n_{\mathrm{it}}$ is the number of iterations and $a$ the number of
fine-level operator actions per iteration. The four times in~\eqref{eq:prob-cost} are recorded per
iteration. The time $t_{\mathrm{op}}$ of one fine-level operator action includes its ghost
exchange, and $t_{\mathrm{pre}}$ covers the rest of the preconditioner (intermediate and coarse
levels, grid transfers and local box solves, with their ghost exchanges). The vector updates
of CG take $t_{\mathrm{vec}}$, and the reductions take $t_{\mathrm{comm}}$: the inner
products of CG and, if present, the coarse-residual sum.

The four times in~\eqref{eq:prob-cost} depend on the implementation and the machine. The first three scale with
$N/P$, while $t_{\mathrm{comm}}$ is bound by latency and grows with $P$. The setup time
$T_{\mathrm{setup}}$ (eigenvalue bounds, coarse operators, one-dimensional eigenproblems) is
reported separately, because it is spent once per mesh and shared by all its solves.

\Cref{tab:results-cost} gives the work and the cost of one iteration. On the octree mesh H-GP on 32
A100 GPUs, CJ-$p$MG needs 257 to 258 times fewer iterations than Jacobi PCG, while each of its
iterations costs 25.1 times more (10.5 against 0.418\,ms). Its CG solver time is therefore 10.2 times
shorter (0.0946 against 0.969\,s). One CJ-$p$MG iteration applies the fine-level operator 13 times,
against once for Jacobi PCG, which explains a factor of 13 in the cost. The intermediate levels,
the transfers and the replicated coarse solve with its global reduction contribute the remaining
factor of about 1.9.
HS-FD needs 11 iterations on this operator, and one of its iterations costs 5.47\,ms on Polaris
and 4.42\,ms on Frontier, 0.52 and 0.62 times a CJ-$p$MG iteration, against 0.418
and 0.319\,ms for Jacobi PCG.

On the uniform Dirichlet mesh U-D on 32 A100 GPUs one FD solve takes 0.785\,ms, the time of 2.6 Jacobi
iterations of 0.299\,ms. Against CJ-$p$MG the advantage of FD splits into a factor of 6 in iterations and a
factor of 5.4 in the cost of one application on Polaris.
The second factor is mostly
algorithmic, since the six one-dimensional mode products of FD cost about as much as two to
three operator actions against 13 for one V-cycle. On the graded mesh G the
iteration factor alone is 250. CJ-$p$MG costs 3.8\,ms per iteration against 1.40\,ms
for one FD solve on Frontier, so nearly all the gain of FD there comes from the iterations.

\subsection{Time to solution}
\label{sec:results-time}

On every tensor-product system the direct FD solve is faster than Jacobi PCG by more than two
orders of magnitude (\cref{tab:results-time}). On the three uniform systems on Frontier it is 16 to 18 times faster than CJ-$p$MG and 11 to 13
times faster than HS-FD, and on Polaris 29 to 32 and 17 to 24 times.
On the graded mesh G on Frontier it is
675 times faster than CJ-$p$MG and 74 times faster than HS-FD.
On Frontier, HS-FD is faster than CJ-$p$MG on all four systems, by factors of 1.32,
1.25, 1.63 and 9.1 on U-D, U-PD, U-P and G. Its cheaper cycle outweighs the one or
two extra iterations it needs on U-D and U-PD.
On 36 Aurora tiles, FD is 38, 14 and $1.8\times10^{3}$ times faster
than CJ-$p$MG and 37, 10 and $2.5\times10^{2}$ times faster than HS-FD on U-PD, U-P and G. 

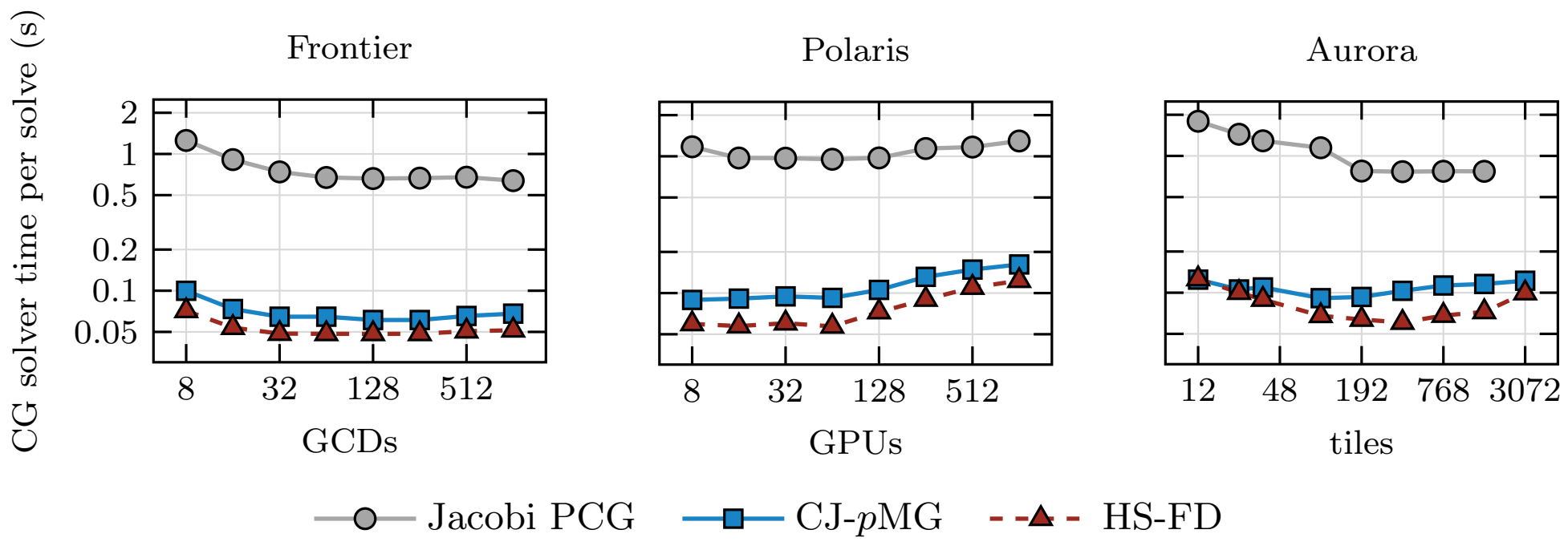
\begin{figure}[tbhp]
\centering
\providecommand{\siscdata}{sisc_figures}
\pgfplotsset{trpaxis/.style={
    sisc/axis, sisc/ranks, sisc/seconds,
    width=0.37\linewidth, height=0.288\linewidth,
    ymin=0.03, ymax=2.5,
    tick label style={font=\footnotesize}, label style={font=\footnotesize}, title style={font=\footnotesize},
  }}
\begin{tikzpicture}
  \begin{loglogaxis}[trpaxis, title={Frontier}, xlabel={GCDs}, ylabel={CG solver time per solve (s)},
      xtick={8,32,128,512}, xticklabels={8,32,128,512},
      legend to name=trplegend, legend columns=3]
    \addplot[sisc/jacobi] table[x=gcds, y=jac_s] {
      gcds,jac_s
      8,1.2550
      16,0.9092
      32,0.7395
      64,0.6730
      128,0.6612
      256,0.6660
      512,0.6763
      1024,0.6380
    };
    \addlegendentry{Jacobi PCG}
    \addplot[sisc/cjpmg]  table[x=gcds, y=cj_s] {
      gcds,cj_s
      8,0.0999
      16,0.0735
      32,0.0646
      64,0.0646
      128,0.0611
      256,0.0611
      512,0.0654
      1024,0.0680
    };
    \addlegendentry{CJ-$p$MG}
    \addplot[sisc/hsfd]   table[x=gcds, y=lv2_s] {
      gcds,lv2_s
      8,0.0713
      16,0.0537
      32,0.0486
      64,0.0484
      128,0.0484
      256,0.0484
      512,0.0506
      1024,0.0514
    };
    \addlegendentry{HS-FD}
  \end{loglogaxis}
\end{tikzpicture}\hfill
\begin{tikzpicture}
  \begin{loglogaxis}[trpaxis, title={Polaris}, xlabel={GPUs}, yticklabels={},
      xtick={8,32,128,512}, xticklabels={8,32,128,512}]
    \addplot[sisc/jacobi] table[x=gpus, y=jac_s] {
      gpus,jac_s
      8,1.1737
      16,0.9720
      32,0.9691
      64,0.9512
      128,0.9735
      256,1.1398
      512,1.1684
      1024,1.2942
    };
    \addplot[sisc/cjpmg]  table[x=gpus, y=cj_s] {
      gpus,cj_s
      8,0.0891
      16,0.0911
      32,0.0946
      64,0.0922
      128,0.1055
      256,0.1314
      512,0.1485
      1024,0.1619
    };
    \addplot[sisc/hsfd]   table[x=gpus, y=lv2_s] {
      gpus,lv2_s
      8,0.0595
      16,0.0574
      32,0.0602
      64,0.0571
      128,0.0727
      256,0.0898
      512,0.1098
      1024,0.1232
    };
  \end{loglogaxis}
\end{tikzpicture}\hfill
\begin{tikzpicture}
  \begin{loglogaxis}[trpaxis, title={Aurora}, xlabel={tiles}, yticklabels={},
      xtick={12,48,192,768,3072}, xticklabels={12,48,192,768,3072}]
    \addplot[sisc/jacobi] table[x=tiles, y=jac_s] {
      tiles,jac_s
      12,1.7933
      24,1.4400
      36,1.2897
      96,1.1472
      192,0.7751
      384,0.7670
      768,0.7724
      1536,0.7718
    };
    \addplot[sisc/cjpmg]  table[x=tiles, y=cj_s] {
      tiles,cj_s
      12,0.1246
      24,0.1056
      36,0.1093
      96,0.0912
      192,0.0931
      384,0.1032
      768,0.1129
      1536,0.1158
      3072,0.1223
    };
    \addplot[sisc/hsfd]   table[x=tiles, y=lv2_s] {
      tiles,lv2_s
      12,0.1259
      24,0.1001
      36,0.0892
      96,0.0678
      192,0.0637
      384,0.0605
      768,0.0681
      1536,0.0724
      3072,0.0994
    };
  \end{loglogaxis}
\end{tikzpicture}

\smallskip
{\small\pgfplotslegendfromname{trplegend}}
\caption[Strong scaling on the octree mesh H-GP]{Strong scaling of the generalized Poisson
solve on the octree mesh H-GP ($N = 11\,390\,280$, hanging nodes,
$\epsilon_\infty = 78.4$, FE degree 7). Frontier: 8 to 1024 MI250X GCDs. Polaris: 8 to 1024 A100
GPUs. Aurora: 12 to 3072 tiles.
CG solver time of one solve with single-reduction CG to a relative true residual of $10^{-10}$, on the
basis of \cref{sec:results-setup}, GPU-aware MPI on Polaris and Frontier, host-staged MPI on
Aurora. CJ-$p$MG needs 9 iterations and HS-FD 11 at every width shown, except 10 for CJ-$p$MG at 36
and 12 for HS-FD at 768 Aurora tiles. Jacobi PCG needs 2317--2321 iterations on Frontier and
2316--2321 on Polaris.}
\label{fig:results-trp-scaling}
\end{figure}

On the octree mesh H-GP, CJ-$p$MG is about 9 to 13 times faster than Jacobi PCG at every
width from 8 to 1024 GCDs, with 9 iterations against 11 for HS-FD (\cref{fig:results-trp-scaling}).
Its CG solver time falls up to 32 GCDs, about $3.6\times10^{5}$ DoFs per GCD, and stays between 0.061
and 0.068\,s up to 1024 GCDs, where about
$1.1\times10^{4}$ DoFs remain per GCD. The V-cycle becomes bound by
latency, so the factor over Jacobi PCG falls from about 13 to about 9 as the width grows. In an instrumented run the wait for ghost data takes 30\% of a V-cycle on 8 GCDs and
57\% on 1024 GCDs, and the replicated coarse solve with its Allreduce 16\% and 29\%.
HS-FD is 1.26 to 1.40 times faster than CJ-$p$MG on Frontier and
1.31 to 1.61 times on Polaris.

The distributed FD solve is faster than the replicated one at every width from 2 ranks on both
machines. On 32 A100 GPUs it takes 0.785 against 9.845\,ms on U-D
(\cref{sec:sm-scaling}), where the strong scaling of the FD solve on U-PD is also given. In the weak scaling on Aurora the
uniform box reaches $N = 1\,573\,037\,747$ DoFs on 3072 tiles, where CJ-$p$MG and HS-FD need
6 and 7 iterations and take 0.786 and 0.852\,s. \Cref{tab:results-aurora-scaling} also gives the
strong scaling on Aurora.

\subsection{Many right-hand sides}
\label{sec:results-mrhs}

Blocks of right-hand sides lower the time per right-hand side of the FD solve by a factor of 3.1
on Frontier and 1.8 on Polaris, and the block FD solve on U-D is 29 to 38 times
faster per right-hand side than block CJ-$p$MG and block HS-FD on Frontier (\cref{tab:results-mrhs}).
For 128 right-hand sides on H-GP on 8 Frontier GCDs, block HS-FD is 1.08 to 1.22
times faster per right-hand side than block CJ-$p$MG for $m = 8$ to 32 (\cref{sec:sm-mrhs}).

\begin{table}[!tbhp]
\caption[Time and iterations on the application problems]{Time in s and, in
parentheses, iterations on the interconnects, the 3D-IC stack, resistivity and EEG, to a relative true residual of $10^{-10}$ at FE degree 4, as in
\cref{sec:results-setup}, on Polaris (P, A100 GPUs), Frontier (F, MI250X GCDs) and Aurora (A, Intel
Max 1550 tiles): CG solver or FD time of one solve, elliptic solve time for nekRS and solve time
without setup for LOR-AMG. Blocks
of $m$ right-hand sides (interconnects, one per conductor, and resistivity): block time divided
by $m$, with the iterations of one right-hand side solved alone. LOR-AMG and nekRS solve one right-hand side at a time and give only these iterations, except nekRS on the interconnects, timed over $m$ sequential solves. Jacobi PCG
iteration ranges are omitted. The preconditioners use FP32 on the fine and intermediate levels and in the HS-FD box
solves. A dash marks a solver not reported.
$^{a}$nekRS multigrid with the coefficient,
Chebyshev--Jacobi smoothing and a Galerkin coarse level, in GMRES(15). $^{b}$HS-FD coarse
correction: BoomerAMG on the low-order-refined operator on F and A, the trilinear level
on P. $^{c}$Without the refinement around the source. $^{d}$nekRS default, preconditioner without the coefficient.}
\label{tab:results-apps}
\centering
\setlength{\tabcolsep}{2pt}
\begin{tabular}{@{}llrrrrrrr@{}}
\toprule
case & & ranks & Jacobi & CJ-$p$MG & HS-FD & FD & LOR-AMG & nekRS \\
\midrule
interconnect & P & 2 & \makecell[r]{0.143\\ (748)} & \makecell[r]{0.0394\\ (18)} & \makecell[r]{0.0191\\ (12)} & --- & (37) & \makecell[r]{0.0768$^{d}$\\ (19)} \\
\quad $m = 4$ & F & 2 & \makecell[r]{0.169\\ (748)} & \makecell[r]{0.0452\\ (18)} & \makecell[r]{0.0255\\ (12)} & --- & (37) & \makecell[r]{0.0883$^{d}$\\ (19)} \\
  & A & 12 & \makecell[r]{0.182\\ (748)} & \makecell[r]{0.055\\ (18)} & \makecell[r]{0.033\\ (12)} & --- & --- & \makecell[r]{0.106$^{a}$\\ (29)} \\
interconnect & P & 64 & \makecell[r]{0.386\\ (1510)} & \makecell[r]{0.151\\ (23)} & \makecell[r]{0.0735\\ (21)} & --- & (37) & \makecell[r]{0.151$^{d}$\\ (19)} \\
\quad $m = 128$\rlap{$^{b}$} & F & 32 & \makecell[r]{0.543\\ (1509)} & \makecell[r]{0.181\\ (23)} & \makecell[r]{0.0750\\ (13)} & --- & (37) & \makecell[r]{0.148$^{d}$\\ (19)} \\
  & A & 24 & \makecell[r]{1.02\\ (1509)} & \makecell[r]{0.497\\ (20)} & \makecell[r]{0.269\\ (13)} & --- & --- & \makecell[r]{0.188$^{d}$\\ (19)} \\
3D-IC & A & 768 & 0.899 & 0.0844 & 0.0732 & 0.0214 & --- & 0.130$^{a}$ \\
\quad $N{\approx}1.0{\times}10^{9}$ & F & 320 & 0.925 & 0.179 & 0.172 & 0.0766 & 0.427 & 0.222$^{a}$ \\
 & P & 320 & 0.959 & 0.355 & 0.335 & 0.0865 & --- & 0.213$^{a}$ \\
\quad $N{\approx}5.0{\times}10^{7}$ & F & 56 & 0.403 & 0.0432 & 0.0379 & 0.00652 & 0.183 & 0.0849$^{a}$ \\
 & A & 48 & 0.633 & 0.0649 & 0.0616 & 0.0126 & --- & 0.106$^{a}$ \\
resistivity & P & 96 & \makecell[r]{9.57\\ (17\,454)} & \makecell[r]{3.08\\ (361)} & \makecell[r]{1.31\\ (43)} & 0.00481 & \makecell[r]{1.23\\ (59)} & \makecell[r]{12.3$^{d}$\\ (756)} \\
\quad $m = 1$ & F & 96 & \makecell[r]{10.1\\ (17\,463)} & \makecell[r]{3.26\\ (361)} & \makecell[r]{0.566\\ (43)} & 0.00618 & \makecell[r]{0.505\\ (60)} & \makecell[r]{9.85$^{d}$\\ (754)} \\
  & A & 96 & \makecell[r]{18.05\\ (17\,456)} & \makecell[r]{18.7\\ (361)} & \makecell[r]{1.77\\ (42)} & 0.0188 & --- & \makecell[r]{12.0$^{d}$\\ (759)} \\
resistivity & P & 96 & 9.20 & 3.27 & 1.07 & 0.00562 & (59) & (756)$^{d}$ \\
\quad $m = 108$ & F & 96 & 10.4 & 3.45 & 0.571 & 0.00611 & (60) & (754)$^{d}$ \\
 & A & 96 & 20.4 & 14.8 & 1.69 & 0.0210 & --- & (759)$^{d}$ \\
EEG, 2\,mm & P & 64 & \makecell[r]{3.92\\ (5341)} & \makecell[r]{0.453\\ (18)} & \makecell[r]{0.370\\ (18)} & --- & \makecell[r]{0.792\\ (38)} & \makecell[r]{0.352$^{a}$\\ (26)} \\
\makecell[l]{EEG, 1\,mm\\ \hspace*{1em}octree} & A & 3072 & \makecell[r]{4.81\\ (10\,775)} & \makecell[r]{1.46\\ (31)} & \makecell[r]{0.725\\ (32)} & --- & --- & --- \\
  & F & 320 & \makecell[r]{9.84\\ (10\,775)} & \makecell[r]{0.644\\ (12)} & \makecell[r]{0.337\\ (17)} & --- & --- & --- \\
\makecell[l]{EEG, 1\,mm\\ \hspace*{1em}uniform$^{c}$} & A & 3072 & \makecell[r]{4.34\\ (10\,243)} & \makecell[r]{1.46\\ (30)} & \makecell[r]{0.654\\ (30)} & --- & --- & \makecell[r]{0.169$^{a}$\\ (23)} \\
  & F & 320 & \makecell[r]{8.87\\ (10\,243)} & \makecell[r]{0.561\\ (11)} & \makecell[r]{0.293\\ (16)} & --- & \makecell[r]{0.489\\ (38)} & \makecell[r]{0.386$^{a}$\\ (22)} \\
\bottomrule
\end{tabular}
\end{table}
\subsection{Application problems}
\label{sec:results-apps}

Application problems extend the tests to further operators, meshes and right-hand sides.
\Cref{tab:results-apps} gives the times and iterations. The
physical setting and validation of each problem are in \cref{sec:sm-apps}.

In capacitance extraction the $k\times k$ bus-crossing interconnect benchmark~\cite{Nabors1991FastCap},
with two crossing layers of $k$ wires, gives one right-hand side per
conductor, a block of $m = 2k$. The conductors are Dirichlet regions inside a graded tensor mesh,
so the operator is not a Kronecker sum. CJ-$p$MG needs 18 to 23 iterations and HS-FD 12 to 21, and
per right-hand side HS-FD is 1.77 and 2.41 times faster than CJ-$p$MG on Frontier and 2.06 and
2.05 times on Polaris with 4 and 128 conductors.

In the thermal analysis of a three-dimensional integrated circuit (3D-IC) the conductivity
$\kappa(z)$ is constant within each layer of the die stack, so the operator remains a Kronecker
sum and the direct FD solve is exact. On the stack of $N = 999\,325\,665$ DoFs it is 3.9 and 3.4 times
faster than CJ-$p$MG and HS-FD on 768 Aurora tiles, about 2.3 times on 320 Frontier GCDs and 4.1 and
3.9 times on 320 Polaris GPUs. On a reduced stack on 56 Frontier GCDs it is 28 times faster than
LOR-AMG.
On 768 Aurora tiles CJ-$p$MG and HS-FD are 11 and 12 times faster than Jacobi PCG, and 9.7 and
10 times on a reduced stack on 48 tiles, where one CJ-$p$MG iteration costs about 17 and
20 Jacobi iterations. Both stacks distribute the coarse level (\cref{sec:precond-common}), and the operator action
sends its larger ghost messages by GPU-aware MPI, while the FD transposes use host memory.

In electrical resistivity tomography (ERT) the three-layer earth of Jaysaval et
al.~\cite{Jaysaval2023PFLOTRAN} is surveyed from 108 electrode positions, a block of 108
right-hand sides on $1.3\times10^{8}$ DoFs. The conductivity depends on $z$ only and the layer interfaces lie on cell faces, so the operator is
a Kronecker sum. The direct FD solve is 82 to 256 times faster than LOR-AMG and 92 to 272 times
faster than HS-FD. On the strongly stretched padding cells CJ-$p$MG needs 361 iterations and HS-FD
42 to 43, and on 96 Frontier GCDs LOR-AMG is 1.12 times faster than HS-FD. Solving the 108
right-hand sides as a block changes the time per right-hand side of CJ-$p$MG and HS-FD by a factor
of 0.79 to 1.06 (\cref{sec:results-ert}).

On the bubbly-flow problem of MacLachlan, Tang and Vuik~\cite{MacLachlan2008BubblyFlow}, a fully
periodic box of 64 spheres with a coefficient ratio of $10^{3}$, HS-FD needs fewer iterations
than CJ-$p$MG and is faster in wall time on all three machines (\cref{tab:sm-bubbly}). The
Taylor--Green vortex and a fiber composite are in \cref{sec:results-tgv,sec:results-fiber}.

\paragraph{Electroencephalography}
\label{sec:results-eeg}

The forward problem of electroencephalography (EEG) gives the scalp potential of a current dipole
in the brain on the New York Head~\cite{HuangNYHead}, with one FE cell per voxel and a
conductivity ratio of about $3.5\times10^{4}$ (\cref{sec:sm-eeg}), so the coefficient does not
separate. On the 1\,mm mesh refined as an octree around the dipole, with hanging nodes ($10^{9}$ DoFs, 3072 Aurora tiles), CJ-$p$MG and HS-FD need 31
and 32 iterations, and HS-FD is 2.01 times faster in wall time, 2.23 times on the uniform 1\,mm
mesh. With the coarse level distributed above $10^{6}$ unknowns, the iterations
do not grow: 18 for both at 2\,mm on 64 Polaris GPUs, and 11
and 16 at 1\,mm on 320 Frontier GCDs, for eight times the DoFs.

\begin{figure}[tbhp]
\centering
\begin{tikzpicture}[font=\footnotesize,every node/.append style={text=black},
    seltest/.style={draw=deepnavy,line width=0.6pt,fill=white,inner sep=3pt,minimum height=0.55cm},
    selleaf/.style={draw=deepnavy,line width=0.6pt,fill=mfblue!35,rounded corners=3pt,
                    inner sep=3pt,minimum height=0.55cm,align=center},
    selfd/.style={selleaf,fill=mforange,line width=1pt},
    seledge/.style={-{Stealth[length=3.5pt]},deepnavy,line width=0.6pt},
    sellab/.style={midway,fill=white,inner sep=1.5pt}]
  \node[seltest] (root) at (5.55,0) {number of right-hand sides};
  \node[seltest] (m1) at (3.0,-1.15) {mesh};
  \node[seltest] (m2) at (9.1,-1.15) {operator a Kronecker sum};
  \node[seltest] (kron) at (1.0,-2.3) {operator a Kronecker sum};
  \node[seltest] (op) at (5.05,-2.3) {operator};
  \node[selfd]   (l1) at (0,-3.55)   {direct FD};
  \node[selleaf] (l2) at (2.0,-3.55) {HS-FD$^{*}$};
  \node[selleaf] (l3) at (4.1,-3.55) {HS-FD or\\ CJ-$p$MG};
  \node[selleaf] (l4) at (6.0,-3.55) {HS-FD};
  \node[selfd]   (l5) at (8.1,-3.55)  {block FD};
  \node[selleaf] (l6) at (10.1,-3.55) {block HS-FD};
  \draw[seledge] (root.south) -- (m1.north) node[sellab] {one};
  \draw[seledge] (root.south) -- (m2.north) node[sellab] {many};
  \draw[seledge] (m1.south) -- (kron.north) node[sellab] {tensor};
  \draw[seledge] (m1.south) -- (op.north) node[sellab] {octree};
  \draw[seledge] (kron.south) -- (l1.north) node[sellab] {yes};
  \draw[seledge] (kron.south) -- (l2.north) node[sellab] {no};
  \draw[seledge] (op.south) -- (l3.north) node[sellab] {$\bK$};
  \draw[seledge] (op.south) -- (l4.north) node[sellab] {$\bK^{\epsilon}$};
  \draw[seledge] (m2.south) -- (l5.north) node[sellab] {yes};
  \draw[seledge] (m2.south) -- (l6.north) node[sellab] {no};
  \node[anchor=north,align=center] at (l1.south) {$\bK$, $\bK + \sigma\bM$,\\ layered $\kappa(z)$};
  \node[anchor=north,align=center] at (l2.south) {conductors};
\end{tikzpicture}
\caption[Selection rule]{Selection rule as a decision tree by the number of right-hand sides, the
mesh class and the operator class. $\bK$: constant coefficient, $\bK^{\epsilon}$: generalized
Poisson operator. Octree: meshes with hanging nodes. The measured margin of each choice is given
in \cref{tab:results-selection}. $^{*}$CJ-$p$MG instead of HS-FD on the battery electrode at a coefficient ratio of $10^{4}$, where
the isolated conducting spheres add near-null modes to the operator.}
\label{fig:results-selection}
\end{figure}
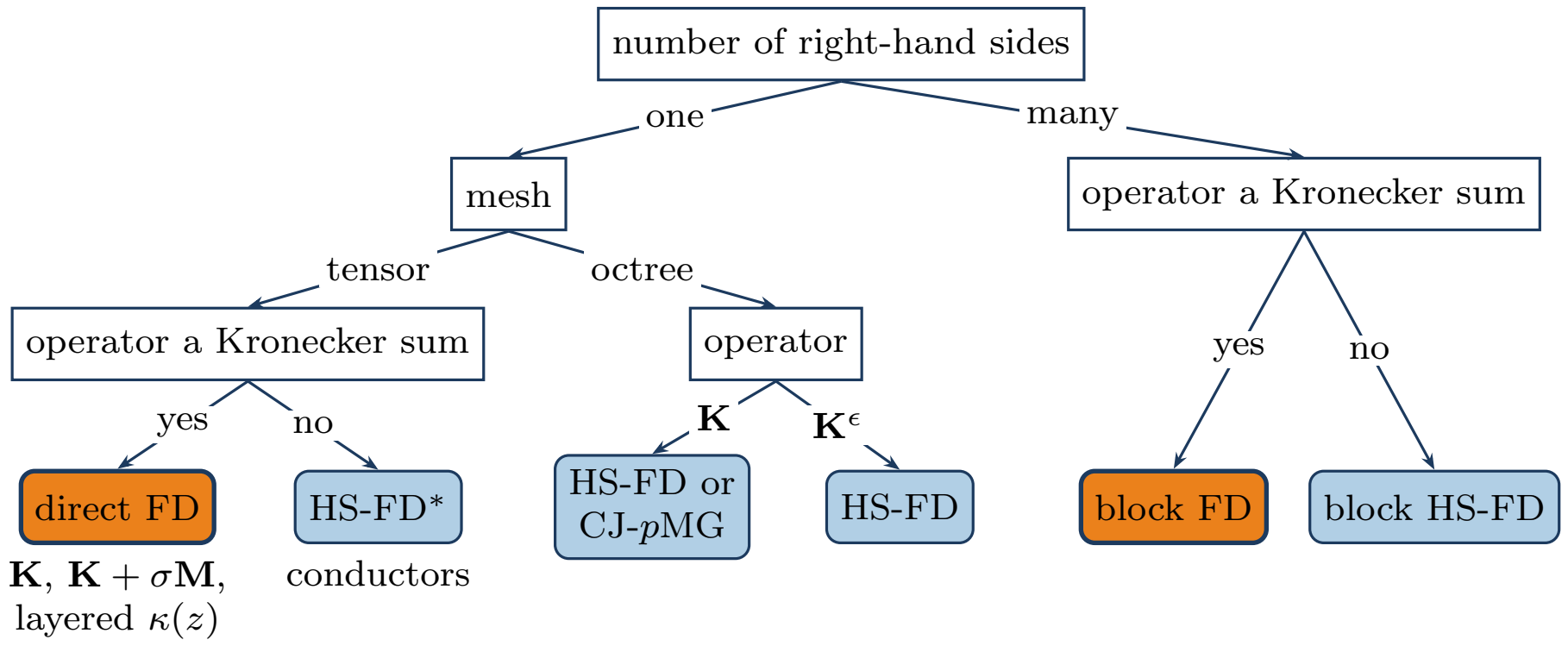
\subsection{Established solvers in other codes}
\label{sec:results-amg}
\label{sec:results-nekrs}

LOR-AMG, CG preconditioned by BoomerAMG~\cite{Henson2002BoomerAMG} on the low-order-refined
operator~\cite{Pazner2020LOR} through MFEM~\cite{Anderson2021MFEM:Library}, solves
the same discrete problems (\cref{sec:sm-comparators}). It needs 31 to 929
iterations, while CJ-$p$MG and HS-FD need 2 to 40 except on stretched cells and on the battery
electrode.

nekRS~\cite{Fischer2022NekRS} is a spectral-element flow solver for GPUs. Its $p$-multigrid
preconditioner for the pressure Poisson problem smooths by Chebyshev-accelerated Schwarz methods
with FD local solves and solves the lowest degree by
BoomerAMG~\cite{Phillips2022TuningSEM}, on the host by default in the released source. On the Poisson box of degree 7 (\cref{sec:sm-comparators}) on 32 A100 GPUs,
CJ-$p$MG and HS-FD are faster than nekRS with the restricted Schwarz smoother of its Kershaw example,
its fastest setting on the box, by factors of 2.8 and 3.6 in the time of their CG iterations and
of 1.8 and 2.2 in the time of their whole solve, from the right-hand side on the device to the
solution on the host, against the elliptic solve time of nekRS.
On 32 Frontier GCDs the factors are 2.2 and 2.9 for the CG iterations and 1.6 and 2.0 for the
whole solve.
On Frontier HS-FD is faster than nekRS with a coefficient-aware multigrid, 0.293 against
0.386\,s on the EEG mesh of $1.0\times10^{9}$ DoFs and 0.172 against 0.222\,s on the 3D-IC stack
(320 GCDs), and 0.0379 against 0.0849\,s on the reduced stack (56 GCDs). On Aurora HS-FD is faster on
the 3D-IC stack, 0.0732 against 0.130\,s (768 tiles), and on the reduced stack, 0.0616 against
0.106\,s (48 tiles). nekRS is faster on the 3D-IC stack on Polaris, 0.213 against 0.335\,s
(320 GPUs), on the EEG mesh on Aurora, 0.169 against 0.654\,s (3072 tiles), and on the 2\,mm EEG
mesh on 64 Polaris GPUs, 0.352 against 0.370\,s (\cref{tab:results-apps}). On the interconnects, which nekRS solves one conductor after another,
HS-FD is faster per right-hand side except with 128 conductors on 24 Aurora tiles, 0.269 against
0.188\,s. On the resistivity mesh the default preconditioner of nekRS, which does not see the
coefficient, needs 754 to 759 iterations against 42 to 43 for HS-FD (12.3 against 1.31\,s on 96
Polaris GPUs).

\subsection{Selection rule}
\label{sec:results-selection}

\Cref{fig:results-selection} condenses the measurements into a selection rule: the direct FD solve
wherever the operator is a Kronecker sum on a tensor-product mesh, and otherwise CJ-$p$MG or HS-FD
by their iterations and the cost of their cycles. HS-FD is selected except on the octree Poisson problem at
$p = 4$, where the two take about the same time, and on the battery electrode at a coefficient
ratio of $10^{4}$, where CJ-$p$MG is 1.47 to 1.84 times faster than HS-FD. On the ERT problem
the direct FD solve is 82 to 256 times faster than LOR-AMG (\cref{sec:sm-selection}).

\section{Conclusions}
\label{sec:discussion}
\label{sec:conclusions}

This paper develops three solvers for high-order elliptic problems in one matrix-free code on
distributed GPUs, written in CUDA, HIP and SYCL. The direct FD solve is distributed, keeps only its
share of the tensor on each MPI rank and also serves as an exact coarse solve. It needs an operator
that is a Kronecker sum, which requires a tensor-product mesh and a coefficient that varies along
one direction only. HS-FD is extended to 2:1 hanging-node octrees with exact cell inverses,
and CJ-$p$MG resolves the hanging-node constraints of every level inside fused kernels. Both
preconditioners use FP32 inside FP64 CG, choose a replicated or a
distributed coarse level, solve blocks of right-hand sides and take every parameter from the
operator.

Where the operator is a Kronecker sum, one FD solve replaces hundreds to hundreds of thousands of
Jacobi PCG iterations, and its time does not depend on the mesh grading up to $R = 4096$. On the
other problems both preconditioners need a number of iterations independent of $h$ under uniform
refinement, and HS-FD keeps its iterations low on strongly graded and octree meshes, with 11 on the
octree generalized Poisson problem on 8 to 1024 MI250X GCDs. CJ-$p$MG and HS-FD are 2.2 to 3.6
times faster than nekRS on the Poisson box in the time of the CG iterations. At $10^{9}$ DoFs nekRS
with a coefficient-aware multigrid is faster than HS-FD on the EEG problem on Aurora and the 3D-IC problem on
Polaris, and slower on the 3D-IC problem on Aurora and Frontier. HS-FD is faster than LOR-AMG on every problem compared except the
ERT problem, where the direct FD solve is faster than both. The measurements give a selection rule
(\cref{fig:results-selection}): the direct FD solve on a Kronecker sum, and otherwise CJ-$p$MG or
HS-FD by iterations and cycle cost.

The iterations of HS-FD grow with the FE degree, as two-level bounds for cell subdomains with
minimal overlap predict~\cite{Toselli2005DomainDecomposition,Brubeck2022VertexStar}, and vertex
patches, which avoid this growth, have no Kronecker form at a hanging interface. Both
preconditioners smooth within one cell, so their iterations grow on stretched cells. For
HS-FD the coarse level of BoomerAMG on the low-order-refined operator corrects this
(\cref{sec:asfd}), for CJ-$p$MG the widened window only in part. Neither coarse space adapts
to the coefficient, so isolated inclusions at a large coefficient ratio still slow
HS-FD.

A replicated coarse level costs one Allreduce per V-cycle
and a coarse solve on every rank, and it grows with $N$, so an AMG coarse level of more than $10^{6}$
unknowns is distributed. The FP32 parts switch to FP64 only at a cell aspect ratio or a ratio of the
coefficient within one cell of 300, not at a coefficient jump between cells. The automatic selection,
timed at setup, can choose differently on other runs or machines and assumes a stationary cycle.

The tensor structure of the direct solve also holds for a separable coefficient $a(x)b(y)c(z)$, for
Robin faces (\cref{sec:fd-kron}) and for the shifted-Laplacian preconditioner of the indefinite
Helmholtz operator, which FD inverts in complex arithmetic. Patches along stretched
cells, coefficient-aware coarse spaces as in black-box multigrid~\cite{MacLachlan2008BubblyFlow}
and the SPE10 permeability field~\cite{Christie2001SPE10} are
future work.

\ifsiscappendix
\appendix
\crefname{appendix}{Appendix}{Appendices}
\Crefname{appendix}{Appendix}{Appendices}
\crefalias{section}{appendix}
\crefalias{subsection}{appendix}
\begingroup
\siscappendixlabels
\siscsupplementbody
\endgroup
\fi

\section*{Code availability}
The solvers are implemented in a branch of DFT-FE~\cite{Motamarri2020DFT-FECalculations,DasGB2023}. The code (branch sisc-iterative-multigrid, commit 5c66ab13e) is archived on Zenodo at \url{https://doi.org/10.5281/zenodo.23797722}. Comparators: MFEM 4.10 (hypre 3.2.0), nekRS v26.0.

\section*{Acknowledgments}
This research used the Polaris and Aurora supercomputers at the Argonne Leadership Computing Facility
(ALCF), a U.S. Department of Energy (DOE) Office of Science user facility at Argonne National
Laboratory, supported by the U.S. DOE Office of Science, Advanced Scientific Computing Research
Program, under Contract No. DE-AC02-06CH11357. This research also used the Frontier
supercomputer at the Oak Ridge Leadership Computing Facility (OLCF), a U.S. DOE Office of
Science user facility at Oak Ridge National Laboratory, supported by the U.S. DOE Office of
Science under Contract No. DE-AC05-00OR22725.

\paragraph{Use of artificial intelligence}
The authors used Claude (Anthropic) to help proofread the manuscript. The authors assume responsibility for all content.

\bibliographystyle{siamplain}
\bibliography{references}

\clearpage
\includepdf[pages=-,noautoscale,offset=0bp -36bp]{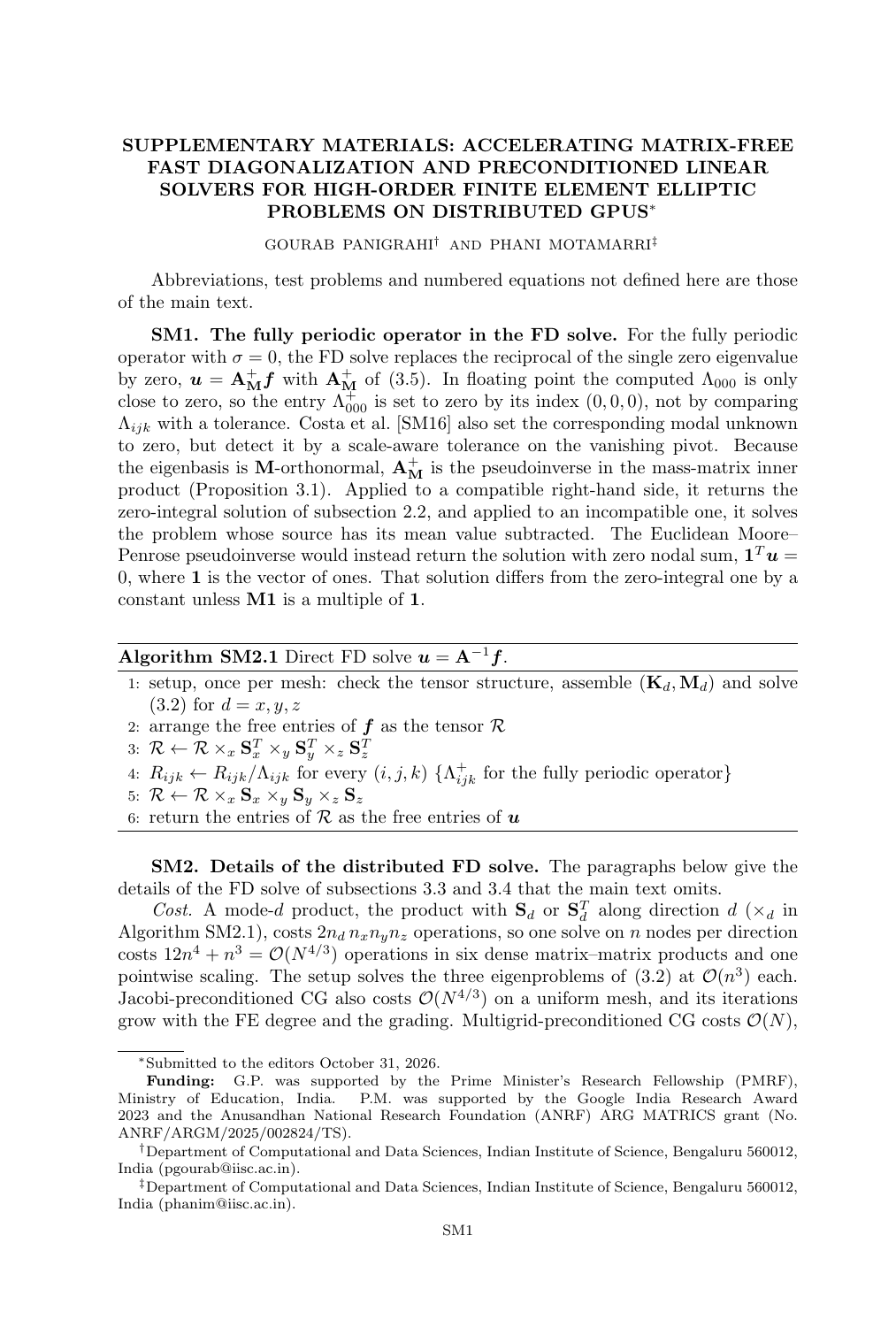}
\end{document}